\documentclass[11pt]{amsart}
\def\Var{{\rm Var}}
\RequirePackage{amssymb}

\newcommand{\bm}[1]{\mbox{\boldmath $#1$}}

\def\1{1\hskip-2.6pt{\rm l}}
\def\N{{\mathbb{N}}}

\def\R{{\mathbb{R}}}

\def\E{{\mathbb{E}}}
\def\P{{\mathbb{P}}}
\def\IL{{{\mathbb{L}}}}
\def\D{{\mathcal D}}

\def\B{{{\mathcal B}}}

\def\F{{\mathcal{F}}}
\def\X{{\mathcal{X}}}
\def\K{{\mathcal{K}}}
\def\T{{\mathcal T}}
\def\I{{\mathcal I}}
\def\J{{\mathcal J}}
\def\M{\mathcal M}
\def\PP{\mathcal P}
\def\EE{{\mathcal E}}
\def\LL{{\mathcal L}}
\def\NN{{\mathcal N}}
\def\RR{{\mathcal R}}
\def\S{{\mathcal S}}
\def\V{{\mathcal V}}
\def\A{{\mathcal A}}

\def\H{{\mathcal H}}
\def\eps{{\varepsilon}}

\newcommand{\pen}{\mathop{\rm pen}\nolimits}
\newcommand{\eref}[1]{(\ref{#1})}
\newcommand{\pa}[1]{\left({#1}\right)}
\newcommand{\cro}[1]{\left[{#1}\right]}
\newcommand{\ab}[1]{\left|{#1}\right|}
\newcommand{\ac}[1]{\left\{{#1}\right\}}
\newcommand{\beqe}{\begin{eqnarray*}}
\newcommand{\eeqe}{\end{eqnarray*}}
\newcommand{\beq}{\begin{eqnarray}}
\newcommand{\eeq}{\end{eqnarray}}
\newcommand{\bq}{\begin{equation}}
\newcommand{\eq}{\end{equation}}
\newcommand{\bin}[2]{\ensuremath{\left(\!\!
\begin{array}{c}
#1\\
#2\\
\end{array}
\!\!\right)}}

\newtheorem{thm}{Theorem}
\newtheorem{lem}{Lemma}
\newtheorem{prop}{Proposition}

\def\PP{\mathcal P}
\def\EE{{\mathcal E}}

\def\S{{\mathcal S}}
\def\A{{\mathcal A}}

\def\I{{\mathcal I}}

\def\H{{\mathcal H}}

\newcommand{\cqfd}{\qquad\framebox[2.7mm]{\rule{0mm}{.7mm}}}

\begin{document}
\title[Estimating the intensity ]{Estimating the intensity of a random measure by histogram 
type estimators} 
\date{April, 2006} 
\author{Yannick Baraud} 
\address{Universit\'e de Nice Sophia-Antipolis, Laboratoire J-A Dieudonn\'e,
  Parc Valrose, 06108 Nice cedex 02} 
\email{baraud@math.unice.fr}
\author{Lucien Birg\'e} 
\address{Universit\'e Paris VI, Laboratoire de Probabilit\'es et 
Mod\`eles Al\'eatoires, bo\^{\i}te 188, 4 place Jussieu, 75252 Paris
Cedex 05} 
\email{lb@ccr.jussieu.fr}
\begin{abstract}
The purpose of this paper is to estimate the intensity of some random measure $N$ on a
set $\X$ by a piecewise constant function on a finite partition of $\X$. Given a (possibly 
large) family $\M$ of candidate partitions, we build a piecewise constant estimator
(histogram) on each of them and then use the data to select one estimator in the
family. Choosing the square of a Hellinger-type distance as our loss function, we show
that each estimator built on a given partition satisfies an analogue of the classical squared
bias plus variance risk bound. Moreover, the selection procedure leads to a final
estimator satisfying some oracle-type inequality, with, as usual, a possible loss
corresponding to the complexity of the family $\M$. When this complexity is not too
high, the selected estimator has a risk bounded, up to a universal constant, by the
smallest risk bound obtained for the estimators in the family. For suitable choices of
the family of partitions, we deduce uniform risk bounds over various classes of
intensities.  Our approach applies to the estimation of the intensity of an inhomogenous Poisson
process, among other counting processes, or the estimation of the mean of a random vector with
nonnegative components.
\end{abstract}
\maketitle
%
\section{Introduction}\label{I}
The aim of the present paper is to design a new model selection procedure in a  statistical 
framework which is general enough to cope simultaneously with the following estimation
problems.

\noindent{\bf Problem 1: Estimating the means of nonnegative data}. 
The statistical problem that initially motivated this research was suggested by Sylvie Huet and 
corresponds to the modeling of data coming from some agricultural experiments. In such an
experiment, the observations are independent nonnegative random variables $N_{i}$ with mean
$s_i$ where $i$ varies among some finite index set $\X$. In this framework, our aim is to
estimate the vector $(s_i)_{i\in\X}$. 

\noindent{\bf Problem 2: Estimating the intensity of a Poisson process}.
We recall that a Poisson process $N$ on the measurable set $(\X,\A)$ with finite mean 
measure $\nu$ is a random measure $N$ on $\X$ such that
\begin{itemize}
\item  for any $A\in\A$, $N(A)$ is a Poisson random variable with parameter $\nu(A)$;
\item for any family $A_1,\ldots,A_n$ of disjoint elements of $\A$, the corresponding
random variables $N(A_1),\ldots,N(A_n)$ are independent. 
\end{itemize}
We can always assume that $\nu$ is finite by suitably  restricting the domain of observation of 
the process. When the mean measure $\nu$ is dominated by some given measure $\lambda$ on
$\X$ then the nonnegative function $s=d\nu/d\lambda$ is called the intensity of $N$. A
Poisson process can be represented as a point process on the set $\X$. Each point represents the
time (if $\X=\R_{+}$) or location of some event. For example, the successive times of failures of
some machine can be represented by a Poisson process on $\X=\R_{+}$. The intensity of the
process models the behaviour of the machine in the following way: the intervals of times on
which the intensity takes large values correspond to periods where failures are expected to be
frequent and in the opposite, those on which the intensity is close to 0 are periods on which
failures are rare. In this statistical framework, our aim is to estimate the intensity $s$ on the basis
of the observation of $N$.

\noindent{\bf Problem 3: Estimating a hazard rate}.
We consider an $n$ sample $T_{1},\ldots,T_{n}$ of non-negative real valued random variables 
with common density $p$ (with respect to the Lebesgue measure on $\R_{+}$) and assume these
to be (possibly) right-censored. This means that there exists i.i.d.\ random variables
$C_1,\ldots,C_n$ such that we actually observe the pairs
$X_{j}=(\widetilde{T}_j,D_j)$ for $j=1,\ldots,n$ with
$\widetilde{T}_j=\min\ac{T_j,C_j}$ and
$D_j=\1_{\{T_j=\widetilde{T}_j\}}$. Such censored data are common in
survival analysis. Typically, $T_{i}$ corresponds to a time of failure or death
which cannot be observed if it exceeds time $C_{i}$. Our aim, here,  is to estimate
the hazard rate $s$ of the $T_{i}$ defined for $t\ge 0$ by $s(t)=p(t)/\P(T_{1}\ge
t)$. 

\noindent{\bf Problem 4: Estimating the intensity of the transition of a Markov process}.
Let $\ac{X_t,\ t\geq 0}$ be a Markov process on $\R_+$ with cadlag paths and a finite number 
of states. We distinguish two particular states, named 0 and 1, and assume that 0 is absorbant and
that there is a positive probability to reach 1. Our aim is to provide an estimation of the intensity
of the transition time $T_{1,0}$ from state 1 to 0. Typical examples arise when 0 means ``death'',
``failure'',\ \dots . An alternative example could be the situation where $T_{1,0}$ measures the
age at which a drug addict makes the transition from soft drugs (state 1) to hard drugs (state 0).
In this case we stop the chain at 0 making this state absorbing.  For $t>0$, we denote by $X_{t-}$
the left-hand limit of the process $X$ at time $t$ and assume that for some measurable
nonnegative function $p$, $\P(T_{1,0}\le t)=\int_{0}^{t}p(u)du$. Note that $p$ is merely the
density of $T_{1,0}$ if  $T_{1,0}<+\infty$ a.s. which we shall not assume. Our aim is to estimate
the transition intensity $s$ of $T_{1,0}$ which is defined for $t>0$ by
$s(t)=p(t)/\P\pa{X_{t-}=1}$. 

For pedagogical reasons mainly, since it has already been extensively studied and can therefore 
serve as a reference, it will be interesting to consider also the much more classical

\noindent{\bf Problem 0: Density estimation}. It is the problem of estimating an unknown density
$s$ from $n$ i.i.d.\ observations $X_1,\ldots,X_n$ with this density. 

All the problems described in the above examples amount to estimating a function $s$ 
mapping $\X$ to $\R_{+}$. For this purpose, we choose a family $\M$ of partitions of $\X$ and
for each $m\in\M$ we design a non-negative estimator $\hat s_{m}$ of $s$ which is constant on
the elements of this partition. We shall call such an estimator an histogram-type estimator. The
performance of $\hat s_{m}$ depends on both $s$ and $m$. Since $s$ is unknown, we cannot
pick the partition which leads to the best estimator. To select a partition in $\M$,
we shall rather use a method solely based on our data leading to some random
partition $\hat m$ and define our resulting estimator as $\hat s_{\hat m}$. Our
objective is to design the selection procedure in such a way that $\hat s_{\hat m}$ performs
almost as well as the best estimator among the family $\ac{\hat s_{m}, m\in \M}$. 

The purpose of this paper is to describe some general setup which allows to 
deal with all the five problems simultaneously, to explain the construction of our
histogram-type estimators $\hat s_{m}$, to design a suitable selection procedure
$\hat m$ and to study the performance of the resulting estimator $\hat s_{\hat
m}$. We shall illustrate our results by numerous examples of family of partitions
and target functions $s$ of interest. For the problems of estimating the intensity
of a Poisson process or a hazard rate on the line, our method provides estimators
than can cope with different families of  functions simultaneously, including
monotone, H\"olderian, or piecewise constant with a few jumps with unknown
locations and sizes. In the multivariate case, we shall also provide some special
method for estimating Poisson intensities with a few spikes with unknown
locations and heights.  

The problem of estimating $s$ by model selection in the first four setups described above did not
receive much attention in the literature with a few noticeable exceptions. Problem~1 is generally
viewed as a regression problem where the mean $s_i$ takes the form $f(x_i)$ for some design
points $x_i$ (typically $f$ is defined on $[0,1]$ and $x_i=i/n$). To perform model selection, one
introduces a wavelet basis and performs a shrinkage of the estimated coefficients of $f$ with
respect to this basis. This amounts to selecting which coefficients will be kept. To this form of
selection pertain the papers by Antoniadis, Besbeas and Sapatinas (2001), Antoniadis and
Sapatinas (2001). Closer to our approach is Kolaczyk and Nowak (2004) based on penalized maximum likelihood. Unlike ours, their approach requires that the means $s_{i}$ be uniformly bounded from above and below by known positive constants. For Problem~2, a similar approach based on wavelet shrinkage is developed in Kolaczyk (1999), but the reference result is Reynaud-Bouret~(2003).
Problems~3 and~4 amount to estimating Aalen's multiplicative intensity $s$ of some counting
process with a bounded number of jumps. The problem of non-parametric estimation of Aalen's
multiplicative intensities has been considered by Antoniadis (1989) who uses penalized
maximum likelihood estimation with a roughness penalty and  gets uniform rates of convergence
over Sobolev balls. Van de Geer (1995) considers the Hellinger loss and establishes uniform
estimation rates for the maximum likelihood estimator over classes of intensities with controlled
bracketting entropy.  Gr\'egoire and Nemb\'e (2000) extend the results of Barron and Cover
(1991) about density estimation to that of intensities. Wu and Wells (2003) and Patil and Wood
(2004) derive asymptotic results for thresholding estimators based on wavelet expansions. All
these results, apart from those of van de Geer, are of an asymptotic nature. Reynaud-Bouret
(2002) introduces a model selection procedure to estimate the intensity. A common feature of
these papers lies in the use of martingales techniques (apart from Gr\'egoire and Nemb\'e,
2000). Unlike theirs, our approach does not require any martingale argument at all.

In Section~\ref{A}, we present a general statistical framework which allows to handle 
simultaneously all the examples we have mentioned. We also make a review of some special
classes of target functions and the various families of models (partitions) to be used in our
estimation procedure. The treatment of our five estimation problems is provided in
Sections~\ref{U} and~\ref{F}. The results presented there derive from a unifying theorem to be
found in Section~\ref{B}. The remainder of the paper is devoted to the most technical proofs. 

In the sequel, we shall make a systematic use of the following notations: constants will be
denoted by $C,C',c,\ldots$ and  may change from line to line; we denote by $\N^*$ the set of
positive integers and we write $x\wedge y$ for $\min\{x,y\}$, $x\vee y$ for $\max\{x,y\}$
and $|m|$ for the cardinality of a set $m$.

\section{Presentation of our method\label{A}}
\subsection{A general statistical framework\label{A1}}
We consider an abstract probability space $(\Omega,\EE,\Bbb{P})$ and a measurable space
$(\X,\A)$ bearing a nonnegative $\sigma$-finite measure $\lambda$. In the sequel $\E$
will denote the expectation with respect to $\Bbb{P}$. We then consider on $\X$ a nonnegative
bounded random process $Y=Y(x,\omega)$, i.e.\ a measurable function from
$\X\times\Omega$ to $\Bbb{R}^+$, and the nonnegative random measure $M$ on $\X$ 
given by $dM=Yd\lambda$.  Besides $M$, we also observe a nonnegative random measure $N$
on $\X$ which satisfies
\begin{equation}
\Bbb{E}[N(A)]=\Bbb{E}\left[\int_A s\,dM\right]<+\infty,\quad\mbox{for all }A\in\A,
\label{Eq-fond}
\end{equation}
for some deterministic nonnegative and measurable function $s$ on $\X$. Note that this 
assumption implies that $N$ is a.s.\ a finite measure. Our aim is to estimate $s$ from the
observations $N$ and $M$. Hereafter, we shall deal with estimators that belong to the
cone $\LL$ of nonnegative measurable functions $t$ on $\X\times\Omega$ such that
$\E\left[\int_{\X}t\,dM\right]<+\infty$. Note that $s$ also belongs to $\LL$. To measure 
the risks of such estimators, we endow $\LL$ with the {\em quasi-distance} (since we may
have $H(t,t')=0$ with $t\ne t'$) $H$  between two elements $t$ and $t'$ of  $\LL$ by 
\[
H^2(t,t')=\int_{\X}\pa{\sqrt{t}-\sqrt{t'}}^2dM,
\]
and set as usual, for $t\in\LL$ and $\F\subset\LL$, $H(t,\F)=\inf_{f\in \F}H(t,f)$.  Given
an  estimator $\hat{s}$ of $s$, i.e.\ a measurable function of $N$ and $Y$ with
$\hat{s}\in\LL$, we define its risk by $\E\left[H^2(\hat{s},s)\right]$.  In most of our
applications, $Y$ is identically equal to 1 in which case $M=\lambda$ is deterministic and if
$t$ and $t'$ are densities with respect to $M$, $H$ is merely the Hellinger distance between
the corresponding probabilities. Only the cases of Problems~3 and~4 require to handle
random measures $M$. 

In order to define our estimators we assume that
\begin{equation}
\Bbb{P}[N(A)>0\mbox{ and }M(A)=0]=0\quad\mbox{for all }A\in\A,
\label{Eq-fond1}
\end{equation}
a property which is automatically fulfilled when $M=\lambda$ is deterministic because 
of~\eref{Eq-fond}. 

\subsection{Histogram-type estimators\label{A2}}
Let us now introduce the histogram-type estimators $\hat s_{m}$ based on some finite 
partition $m$ of $\X$. We consider the subset $\J=\{A\in\A\,|\,\E[M(A)]<+\infty\}$ of
$\A$ and define the {\em model} $S_m$ as the set of (possibly random) nonnegative piecewise constant functions on $\X$:
\[ 
S_m=\left\{\left.t=\sum_{I\in m\cap \J}t_I\1_{I}\,\right|\,
t_I=t_I(\omega)\in\R\mbox{ for all }I\in m, \omega\in\Omega\right\}\bigcap\LL.
\] 
We then define the histogram estimator $\hat{s}_m$ as the element of  $S_m$ given (with
the convention $0/0=0$) by 
\[
\hat{s}_m=\sum_{I\in m\cap \J}\frac{N(I)}{M(I)}\1_{I}.
\]
Note that $\hat{s}_m$ is a.s.\ well-defined because of (\ref{Eq-fond1}). We shall, hereafter, call
it the {\em histogram estimator} based on $m$. 

Under suitable assumptions that will be satisfied for Problems~0, 1 and~2 (the case of hazard
rates and Markov processes being more complicated), we shall prove for $\hat{s}_m$ a risk bound
of the form
\begin{equation}
\E\cro{H^2(\hat{s}_m,s)}\le C_0\left\{\E\cro{\pa{H^2(s,S_{m})}}+C_P|m|\right\},
\label{Eq-risk9}
\end{equation}
where $C_0$ is a numerical constant and $C_P$ depends on the problem we consider. For instance, $C_P=n^{-1}$ for density estimation and $C_P=1$ for estimating the intensity of a Poisson process. We recover here the usual decomposition of the risk bounds into
an approximation term which involves the distance of the parameter from the model and a
complexity term proportional to the number $|m|$ of parameters that describe the model.
\subsection{The selection procedure}\label{A3}
Given the family of models $\{S_m,m\in\M\}$ corresponding to a finite or countable family 
$\M$ of partitions $m$, we consider, in order to define our model selection procedure, the
possibly enlarged family 
\[
\overline{\M}=\{m\vee m'\mbox{ for }m,m'\in\M\};\quad
m\vee m'=\{I\cap I'\,|\,I\in m, I'\in m',\ I\cap I'\neq\varnothing \},
\]
so that $m\vee m'$ is again a finite partition of $\X$.

We shall systematically make the following assumption about the family $\M$.\vspace{2mm}\\
{\bf H}\,:\, {\it There exists some $\delta\ge1$ such that
$|m\vee m'|\le\delta \pa{|m|+|m'|}$ for all $(m,m')\in\M^2$.}\vspace{2mm}\\
We then introduce a penalty function ``pen" from $\M$ to $\R_+$ to be
 described below and,
for $m\ne m'\in\M$ we consider the test statistic
\begin{equation}
T_{m,m'}(N)=H^2(\hat{s}_m,\hat{s}_{m\vee m'})-H^2(\hat{s}_{m'},\hat{s}_{m\vee m'})+
16[\pen(m)-\pen(m')].
\label{Eq-def2}
\end{equation}
The corresponding test between $m$
and $m'$ decides $m$ if $T_{m,m'}<0$, $m'$ if $T_{m,m'}>0$ and at random if
$T_{m,m'}=0$. Note that the tests corresponding to $T_{m,m'}$ and $T_{m',m}$
are the same. We then set, for all $m\in\M$,
\[
\RR_m=\{m'\in\M,\;m'\ne m\,|\,\mbox{the test based on }T_{m,m'}\mbox{ rejects }m\}
\]
and, given some $\varepsilon>0$, we define $\hat{m}$ to be any point in $\M$ such that
\begin{equation}
\D(\hat{m})\le\inf_{m\in\M}\D(m)+\varepsilon/3\qquad\mbox{with}\qquad
\D(m)=\sup_{m'\in\RR_m}\left\{H^2(\hat{s}_m,\hat{s}_{m'})\right\}.
\label{Eq-def1}
\end{equation}
This model selection procedure results in an estimator $\tilde{s}=\hat{s}_{\hat{m}}$ that
we shall call {\em penalized histogram estimator} (in the sequel PHE, for short) based on
the family of models $\{S_m,m\in\M\}$ and the penalty function $\pen(\cdot)$. As to the
penalty, it is the sum of two components: $\pen(m)=c_1|m|+c_2\Delta_m$ with $c_1$ and
$c_2$ depending on the framework and $\Delta_m$ being a nonnegative weight associated to
the model $S_m$. We require that those weights satisfy
\begin{equation}
\sum_{m\in\M}\exp[-\Delta_m]=\Sigma<+\infty.
\label{som}
\end{equation}
If $\Sigma=1$, the choice of the $\Delta_m$ can be viewed as the choice of a prior distribution
on the models. For related conditions and their interpretation, see Barron and Cover (1991),
Barron, Birg\'e and Massart (1999) or Birg\'e and Massart (2001). The constant 16
in~\eref{Eq-def2} plays no particular role and has only been chosen in order to improve the
legibility of our main results. Our selection procedure can be viewed as a mixture between a
method due to Birg\'e (1983 and 2006) based on testing and an improved version of the original
Lepski's method, as described in Lepski (1991) and subsequent work of the same author. This
improved version was presented by Lepski in a series of lectures he gave at Garchy in 1998.
\subsection{Risk bounds for the procedure}\label{A4}
As we shall see later, with a suitable choice of $\varepsilon$, the performances of this procedure
for Problems~0, 1 and 2 are described by risk bounds of the following form:
\begin{equation}
\E\cro{H^2(\tilde{s},s)}\le C'_0\inf_{m\in\M}\left\{\E\cro{\pa{H^2(s,S_{m})}}+
C_P|m|\left[1+|m|^{-1}\left(\Delta_m+\Sigma^2\right)\right]\right\},
\label{Eq-risk8}
\end{equation}
where $C'_0$ is a numerical constants and $C_P$ as in (\ref{Eq-risk9}). Comparing 
(\ref{Eq-risk8}) with (\ref{Eq-risk9}), we see that the estimator $\tilde{s}$ achieves a risk bound
comparable, up to a constant factor, with the best risk bound obtained by the estimators
$\hat{s}_m$ provided that $\Sigma$ is not large and
$\Delta_m$ not much larger than $|m|$. Note that these two restrictions are, to some extent,
contradictory since the smaller $\Delta_m$, the larger $\Sigma$, although it is clearly
unnecessary to choose $\Delta_m$ smaller than $|m|$. Therefore, if
$\sum_{m\in\M}e^{-|m|}$ is not large, one can merely take $\Delta_m=|m|$. Otherwise, the
choice of the $\Delta_m$ will be more delicate but we should keep in mind that, if $\Sigma$ is
not large, the performance of $\tilde{s}$ will be as good (up to a constant factor) as the
performance of any $\hat{s}_m$ for which $\Delta_m\le|m|$.

\section{A review of the models we shall use\label{V}}
\subsection{Some classes of functions of special interest\label{V1}}
The motivations for the choice of some family of models $\{S_m,m\in\M\}$ are twofold. First,
there is the restriction that $\M$ should satisfy Assumption~{\bf H} and there are two main
examples of such families. In the "nested" case, the family is totally ordered for the inclusion and
thus, we either have $m\vee m'= m$ or $m\vee m'=m'$ for all $m$ and  $m'$ in $\M$. Then,
$\overline{\M}=\M$ and $\delta=1$. Another situation where Assumption~{\bf H} is
satisfied with $\delta=1$ occurs when $\X$ is either $\R$ or some subinterval of $\R$ and each
$m\in\M$ is a finite partition of $\X$ into intervals.

The second motivation is connected to the approximation properties of the models. If, for
instance, we believe that the true $s$ is smooth or monotone, one should introduce families of
models that approximate reasonably well such functions. In the sequel, we shall put a
special emphasis on the following classes of functions:
\begin{itemize}
  \item {\it Monotone functions.} For $\X$ an interval of $\R$ with interior $\overset{\circ}{\X}$ and $R$ a
positive number, we denote by $\S^1(R)$ the set of monotone functions $t$ on $\X$ such that
$\sup_{x,y\in\overset{\circ}{\X}}|t(x)-t(y)|\le R$.
  \item {\it Continuous functions.} Let $w$ be a modulus of continuity on $[0,1)$, i.e.\ a
continuous nondecreasing function with $w(0)=0$ --- see additional details in DeVore and
Lorentz (1993) ---. We denote by ${\S}^2(w)$ the set of functions $t$ on $[0,1)$ such that
$|t(x+y)-t(x)|\le w(y)$ for all $x\in[0,1)$ and $0\le y\le 1-x$. For $0<\alpha\le1$ and
$R>0$, the H\"older class $\H^R_\alpha$ is the class ${\S}^2(w)$ with $w(y)=Ry^\alpha$.
More generally we say that a function $u$ defined on $\V\subset[0,1)^k$ for some $k\ge1$
belongs to the set $\H^R_{\alpha}(\V)$, $\alpha\in]0,1)$, $R>0$, if 
\[ 
|u(x)-u(y)|\le R\sum_{j=1}^k|x_j-y_j|^\alpha\quad\mbox{for all } x,y\in\V.
\] 
  \item {\it Piecewise constant functions.} If the function $t$ defined on $[0,1)$ is constant over
some intervals and then jumps from time to time, it is a piecewise constant function of the form
\begin{equation}
t=\sum_{k=1}^D t_k\1_{[x_{k-1},x_k)}\quad\mbox{with }0=x_0<x_1<\ldots<x_D=1.
\label{eq:s}
\end{equation}
We shall denote by $\S^3(D,R)$ the class of such piecewise functions such that
$\sup_{1\le k\le D}t_k\le R$. Note that this would correspond to a parametric
model with $D$ parameters if the locations of the jumps were known. We shall
restrict our attention to $D\ge2$ since $\S^3(1,R)$ only contains constant
functions and is then a subset of $\S^2(w)$ with $w\equiv0$.
  \item {\it Besov balls and functions of bounded variation.}
Here we consider functions $t$ defined on $[0,1)$. Given positive numbers $\alpha, p$ and
$R$, we denote by $\B^\alpha_{p,\infty}(R)$, the closed Besov ball of radius $R$ centered
at zero of the Besov space $B^\alpha_{p,\infty}([0,1))$, i.e.\ the set of functions $t$ in this
space with Besov semi-norm $|t|_{B^\alpha_{p,\infty}}\le R$. Analogously, we set
$\B_{BV}(R)$ for the set of functions $t$ of bounded variation with $\Var^*(t)\le R$.
We refer to Chapter~2 of the book by DeVore and Lorentz (1993) for details on Besov spaces
and the definition of Besov semi-norms, functions of bounded variation and the
variation semi-norm $\Var^*$. Note that $\S_1(R)\subset \B_{BV}(R)$. We shall also consider
the multidimensional Besov spaces $B^\alpha_{p,\infty}([0,1)^k)$ for $k\ge2$.
\end{itemize}
\subsection{Some typical models\label{V2}}
Let us now describe a few useful families of models and corresponding choices for the weights
$\Delta_m$ that satisfy (\ref{som}).

\subsubsection{Example 1: models for functions on $[0,1)$}\label{A4a}
The following models are suitable for approximating functions belonging to the classes that we
just mentioned. Since they are based on partitions of $[0,1)$ into intervals, they satisfy
Assumption~{\bf H} with $\delta=1$. Let ${\J}_l=\left\{j2^{-l}, j\in\N\right\}$ and
$\J_\infty=\cup_{l\in\N}{\J}_l$ be the set of all dyadic points in $[0,1)$. To build $\M$, we
consider partitions $m=\{I_1,\ldots,I_D\}$ of $[0,1))$ generated by increasing sequences
$\{0=x_0<x_1<\ldots<x_D=1\}$ with $I_i=[x_{i-1},x_i)$. We  then define $\M$ to be the set of
all such partitions with $x_i\in\J_\infty$ for $1\le i\le D-1$. Therefore, whatever
$m\in{\M}$, the elements of $S_m$ are piecewise constant functions  with $D$ pieces and
jumps located on the grid ${\J}_\infty$. The novelty of this particular family of  partitions lies
in the fact that there is no lower bound on the length of the intervals on which the partitions
are built. It will be useful to single out the set $\M_R=\{m_k,k\in\Bbb{N}\}$ of regular
dyadic partitions where $m_k$ is the partition of $[0,1)$ into $2^k$ intervals of length
$2^{-k}$. In particular, $m_0=[0,1)$.

One possible way of defining the corresponding weights $\Delta_m$ is as follows. For
$l\in\N^\star$ and $2\le D\le2^l$ we define ${\M}_{l,D}$ as the set of all partitions
$m$ with $|m|=D$ and $l$ is the smallest integer such that
$\{x_1,\ldots,x_{D-1}\}\subset{\J}_l$. Then, ${\M}=
\left[\bigcup_{l\ge1}\left(\bigcup_{D=2}^{2^l}{\M}_{l,D}\right)\right]\bigcup\{m_0\}$.
We choose $\Delta_{m_0}=1$ and
\begin{equation}
\Delta_m=D(l\log 2+2-\log D)+2\log l\quad\mbox{if }m\in{\M}_{l,D}.
\label{Eq-Dm4}
\end{equation}
Since $|{\M}_{l,D}|\le\binom{2^l-1}{D-1}\le\binom{2^l}{D}\le(2^le/D)^D$, we derive
from (\ref{Eq-Dm4}) that
\begin{eqnarray*}
\sum_{m\in\M\setminus\{m_0\}}\exp[-\Delta_m]&<&\sum_{l\ge1}\;\sum_{D=2}^{2^l}
|{\M}_{l,D}|l^{-2}\exp[-D(l\log 2+2-\log D)]\\&\le&\sum_{l\ge1}\;
\sum_{D\ge2}l^{-2}e^{-D}\;\;=\;\;\frac{\pi^2-6}{6e(e-1)}\;\;<\;\;0.14
\end{eqnarray*}
and it follows that (\ref{som}) is satisfied. 
\subsubsection{Special partitions derived from adaptive approximation algorithms}\label{A4w}
It is easily seen that the family $\M$ of partitions we introduced for Example~1 is too rich
for choosing $\Delta_m=c|m|$ for all $m$ and $c$ a fixed constant since then (\ref{som})
would not be satisfied. For partitions in $\M_{l,D}$ with $l>D$, $\Delta_m$ behaves as $l|m|$
and $l$ can be arbitrarily large. Fortunately, there exists a subset $\M^1_T$ of $\M$, which is of
special interest because of its approximation properties with respect to functions in Besov
spaces, and such as it is possible to choose $\Delta_m=2|m|$ for $m\in\M^1_T$. This will
definitely improve the performances of the PHE for estimating functions in Besov spaces. Let
us now describe $\M^1_T$.

Among all partitions on $[0,1)$ with dyadic endpoints, some of them, which are in one-to-one
correspondance with the family of complete binary trees, can be derived by the following
algorithm described in Section~3.3 of DeVore (1998). One starts with the root of the tree which
corresponds to the interval $[0,1)$ and decides to divide it into two intervals of length $1/2$ or
not. We assume here that all intervals contain their left endpoint but not the right one. If one
does not divide, the algorithm stops and the tree is reduced to its root. If one divides, one gets
two intervals corresponding to adding two sons to the root. Then one repeats the procedure with
each interval and so on\dots . At each step, the terminal nodes of the tree correspond to the
intervals in the partition and one decides to divide any such interval into two equal parts or not.
Dividing means adding two sons to the corresponding terminal node. The whole procedure
stops at some stage producing a complete binary tree with $D$ terminal nodes and the
corresponding partition of $[0,1)$ into $D$ intervals. This is the type of tree which comes out of an algorithm like CART, as described by  Breiman et
al.\ (1984). Such constructions and the corresponding selection procedure resulting from the
CART algorithm have been studied by Gey and Nedelec (2005).  We denote by $\M^1_T$ the
subset of $\M$ of all partitions that can be obtained in this way. Note here that the set $\M_R$ of
regular partitions is a subset of $\M^1_T$.

It is known that the number of complete binary trees with $j+1$ terminal nodes is given by the
so-called Catalan numbers $(1+j)^{-1}\bin{2j}{j}$ as explained for instance in Stanley (1999,
page 172). As a consequence, we can redefine $\Delta_m=2|m|$ for $m\in\M^1_T$ and, using
the fact (which derives from Stirling's expansion) that $\bin{2j}{j}\le4^j$, get
\begin{eqnarray*}
\sum_{m\in\M^1_T}\exp[-\Delta_m]&<&\sum_{j\ge0}\sum_{\{m\in\M^1_T\,|\,|m|=
1+j\}}\exp[-2(j+1)]\\&=&\sum_{j\ge0}\frac{\bin{2j}{j}\exp[-2(j+1)]}{j+1}\;\;\le\;\;e^{-2}
\sum_{j\ge0}\frac{(2/e)^{2j}}{j+1}\;\;=\;\;\Sigma'_1.
\end{eqnarray*}
Finally (\ref{som}) is satisfied with $\Sigma<\Sigma'_1+0.14$.

\subsubsection{Example 2: estimating functions with radial symmetry\label{A4b}}
There are situations where one may assume that the value of $s(x)$ only depends on the 
Euclidean distance $\|x\|$ between this point and some origin in which case one can write 
$s(x)=\Phi(\|x\|)$. In such a case, it is natural to estimate $s$ on a ball, which we may assume,
without loss of generality, to be the open unit ball $\B_k$ of $\R^k$. To any partition $m$ of
$[0,1)$ we can associate a partition of $\B_k$ with elements $J=\{x\,|\,\|x\|\in I\}$ where $I$
denotes an element of $m$. For simplicity, we shall identify the two partitions (the first one
of $[0,1)$ and the new one of $\B_k$) and denote both of  them by $m$. In the sequel, we
shall focus our attention on the family of partitions of Example~1 with the weights defined in 
Section~\ref{A4w}.

\subsubsection{Example 3: estimating functions on $[0,1)^k, k\ge2$}\label{A4c}
To deal with the case $\X=[0,1)^k$, let us first introduce some notations. For $j\in\N$ we
consider the set
\[
\NN_j=\left\{\left.\bm{l}=(l_1,\ldots,l_k)\in\N^k\,\right|\,1\le l_i\le
2^j\;\;\mbox{for}\;\;1\le i\le k\right\}
\]
and for $j\in\N$ and $\bm{l}\in\NN_j$ the cube $K_{j,\bm{l}}$ given by
\[
K_{j,\bm{l}}=\left\{\left.\bm{x}=(x_1,\ldots,x_k)\in[0,1)^k\,\right|\,(l_i-1)2^{-j}\le
x_i<l_i2^{-j}\;\;\mbox{for}\;\;1\le i\le k\right\}.
\]
We set $\K_j=\left\{K_{j,\bm{l}},\bm{l}\in\NN_j\right\}$ and $\K=\bigcup_{j\ge0}\K_j$. 

Let $\PP$ be the collection of all finite subsets $p$ of $\K\setminus\K_0$ consisting of disjoint
cubes. To each $p\in\PP$, we associate the positive quantity $J(p)=\inf\{j\,|\,p\cap\K_j\ne
\varnothing\}$ ($J(\varnothing)=+\infty$) and the partition $m_p$ generated by $p$, i.e.\ $m_p=
\ac{I\in p}\bigcup\ac{[0,1)^{k}\setminus \cup_{I\in p}I}$ provided that this last set is not
empty and $m_p=\ac{I\in p}$ otherwise. We finally set $\M=\{m_p\vee\K_j\mbox{ with }
p\in\PP\mbox{ and }j<J(p)\}$. Note here that the mapping $(j,p)\mapsto m_p\vee\K_j$ is not
one to one. For instance $m_\varnothing\vee\K_j=\K_j=\K_j\vee\K_{j-1}$. We shall prove in
Section~\ref{P9} the following result:
%
\begin{lem}\label{L-Ex4}
The family $\M$ satisfies Assumption~{\bf H} with $\delta=2$.
\end{lem}
In order to define the weights $\Delta_m$, we shall distinguish a special subset $\M^k_T$ of
$\M$ which is the $k$-dimensional analogue  of the one we considered in Section~\ref{A4w}.
Here one starts the algorithm with $\X=[0,1)^k$ (which corresponds to the root of the tree) and
at each step get a partition of $\X$ into a finite family of disjoint cubes of the form
$K_{j,\bm{l}}$. One then decides to divide any such cube into the $2^k$ elements of
$\K_{j+1}$ which are contained in it or not. Again, this corresponds to growing a complete
$2^k$-ary tree, partioning a cube meaning adding $2^k$ sons to a terminal node and the set 
$\M^k_T$ of all partitions that can be constructed in this way corresponds to the set of complete
$2^k$-ary trees. As for $k=1$, $\M^k_T$ contains the set
$\M_R=\{m_\varnothing\vee\K_j,\ j\ge 0\}$ of all regular partitions of $\X$ into $2^{kj}$
cubes of equal volume. Working with
$\M$ instead of the much simpler family $\M_R$ allows to handle less regular functions like
those which have a few spikes or are less smooth on some subset of $\X$.

If $m\in\M^k_T$ we take $\Delta_m=|m|$ and otherwise we set
\[ 
\Delta'_{j,p}=j+k\sum_{i\ge1}(j+i)\ |p\cap\K_{j+i}|\quad\mbox{ for }p\in\PP
\mbox{ and }j<J(p)
\]
and 
\begin{equation}
\Delta_m=\inf_{\{(j,p)\,|\,m=m_p\vee\K_j\}}\ac{\Delta'_{j,p}}
\quad\mbox{for }m\in\M\setminus\M^k_T.
\label{Eq-delta'}
\end{equation}
Note that the ratio $\Delta_m/|m|$ is unbounded for $m\not\in\M^k_T$ as shown by the
example of $m=m_p\vee\K_0$ with $p$ reduced to a single element of $\K_j$, $j>0$. Then
$|m|=2$ while $\Delta_m=kj$ may be arbitrarily large. For the partitions $m$  belonging to
$\M^k_T$ we use the fact --- see Stanley (1999) --- that any complete $l$-ary tree has a number of
terminal nodes of the form $1+j(l-1)$ for some $j\in\Bbb{N}$ and that the number of such trees
with $1+j(l-1)$ terminal nodes is $[1+j(l-1)]^{-1}\bin{lj}{j}$. For $l=2^k$ we derive that the
number of partitions in $\M^k_T$ with $1+j(2^k-1)$ elements is $[1+j(2^k-1)]^{-1}\bin{2^kj}{j}$.
Moreover, since $k\ge2$, we check that
\[
\Delta_m>j(k\log2+1)+\log(j+1)\quad\mbox{if }|m|=1+j(2^k-1).
\]
Since $\bin{lj}{j}\le(le)^j$, it follows that
\begin{eqnarray*}
\sum_{m\in\M^k_T}\exp[-\Delta_m]&<&\sum_{j\ge0}\sum_{\{m\in\M^k_T\,|\,|m|=
1+j(2^k-1)\}}
\frac{\exp[-j(k\log2+1)]}{j+1}\\&=&\sum_{j\ge0}\frac{\bin{2^kj}{j}\left(2^ke\right)^{-j}}
{(j+1)[1+j(2^k-1)]}\\&\le&\sum_{j\ge0}\frac{1}{(j+1)[1+j(2^k-1)]}\;\;=\;\;\Sigma'_k.
\end{eqnarray*}
Let us now turn to the partitions of the form $m_p\vee\K_j$. For such a partition $p\cap\K_{j'}
=\varnothing$ for $j'\le j$ and, for $i\ge1$, $|p\cap\K_{j+i}|=l_i$ with $0\le
l_i\le2^{k(j+i)}$. Moreover, the number of those $p\in\PP$ such that $|p\cap\K_{j+i}|=l_i$
for a given sequence
$\bm{l}=(l_i)_{i\ge1}$ with a finite number of nonzero coefficients is bounded by 
$\prod_{i\ge1}\bin{2^{k(j+i)}}{l_i}$. It follows from (\ref{Eq-delta'}) that
\begin{eqnarray*}
\sum_{m\in\M'}\exp[-\Delta_m]&\le&\sum_{j\ge0}\sum_{\{p\in\PP\,|\,J(p)>j\}}e^{-j}
\prod_{i\ge1}e^{-k(j+i)|p\cap\K_{j+i}|}\\&\le&\sum_{j\ge0}e^{-j}\sum_{\bm{l}}\,
\sum_{\{p\,|\,|p\cap\K_{j+i}|=l_i\mbox{ for }i\ge1\}}\,\prod_{i\ge1}e^{-k(j+i)l_i}\\&\le
&\sum_{j\ge0}e^{-j}\sum_{\bm{l}}\prod_{i\ge1}\bin{2^{k(j+i)}}{l_i}e^{-k(j+i)l_i}\\&\le&
\sum_{j\ge0}e^{-j}\prod_{i\ge1}\sum_{l_i=0}^{2^{k(j+i)}}\bin{2^{k(j+i)}}{l_i}e^{-k(j+i)l_i}\\&
=&\sum_{j\ge0}e^{-j}\prod_{i\ge1}\left(1+e^{-k(j+i)}\right)^{2^{k(j+i)}}\\&=&\sum_{j\ge0}
\exp\left[-j+\sum_{i\ge1}2^{k(j+i)}\log\left(1+e^{-k(j+i)}\right)\right]\\&\le&\sum_{j\ge0}
\exp\left[-j+\sum_{i\ge1}(e/2)^{-k(j+i)}\right]\;\;=\;\;\Sigma''_k\;\;<\;\;+\infty.
\end{eqnarray*}
Finally  we can conclude that (\ref{som}) holds with $\Sigma<\Sigma'_k+\Sigma''_k$. 

\subsubsection{Models for $n$-dimensional vectors}\label{A4z}
To handle the problem we started with in the introduction, we may assume that our finite
index set $\I$ is actually  $\X=\{1,\ldots,n\}$, the estimation of the function $s$ from
$\X$ to $\R_+$ amounting to the estimation of the vector $(s_1,\ldots,s_n)^t\in\R_+^n$ with
coordinates $s_i=s(i)$. 

\subsubsection*{Example 4}
If one assumes that either $s_i$ varies smoothly with $i$ or is monotone or piecewise
constant with a small number of jumps, it is natural to choose for $m$ a partition of $\X$ into
intervals and for $\M$ the set of all such partitions. Note that this family satisfies
Assumption~{\bf H} with $\delta=1$. Setting here $\Delta_m=|m|+\log\bin{n-1}{|m|-1}$,
we get (\ref{som}) with $\Sigma<(e-1)^{-1}$ since there are $\bin{n-1}{D-1}$ partitions in $\M$
with $D$ elements for $1\le D\le n$. 

\subsubsection*{Example 5}
An alternative case is the case when $s$ is constant, equal to $\overline{s}$ on $\X$ except for
a few number of locations $i$ where $s(i)\ne\overline{s}$. Since the number $k$ of such
locations is unknown, it is natural, for each $k\in\{0,\ldots,n-1\}$ to define $\M_k$ as the set
of partitions of $\X$ with $k$ singletons and the set of the $n-k$ remaining points. We finally
set $\M=\cup_{0\le k\le n-1}\M_k$. Then Assumption~{\bf H} holds with $\delta=1$. For 
$m\in\M_k$, $|m|=k+1$ and we set $\Delta_m=\log\bin{n}{k}+k=
\log\bin{n}{|m|-1}+|m|-1$, so that (\ref{som}) holds with $\Sigma<e/(e-1)$.

\section{The case of a deterministic measure $M$}\label{U}
Let us now see how our general framework applies to Problems~1 and~2. Besides these, our
setup also covers the problem of density estimation. Although there is a huge amount of
literature on density estimation, our method brings some improvements to known results on
partition selection for histograms. Moreover, since this problem has attracted so much
attention, it can serve as pedagogical example and reference for the sequel. This is why,
before considering more original and less studied frameworks, we shall start our review by
this quite familiar estimation problem.

\subsection{Density estimation}\label{U2}
We consider the classical problem of estimating an unknown density $s$ from a
sample of size $n$, which means that we have at hand an i.i.d.\ sample $X_1,\ldots,X_n$
from a distribution with unknown density $s$ with respect to some given measure
$M=\lambda$ on $\X$. We define $N$ to be the empirical distribution: $N(A)=n^{-1}
\sum_{i=1}^n\1_{X_i\in A}$. Then, as required, $\E\cro{N(A)}=\int_A s\,d\lambda$ for all
measurable subsets $A$ of $\X$. In this case the distance $H$ is merely a version of the
Hellinger distance between densities.

Within this framework, we can prove the following general result.
%
\begin{thm}\label{densite}
Assume that the family $\M$ satisfies Assumption~{\bf H} and the weights $\ac{\Delta_m,\
m\in\M}$ are chosen so that \eref{som} holds. Then the penalized histogram estimator
$\tilde{s}=\hat s_{\hat m}$ defined in Section~\ref{A3} with $\pen(m)\ge n^{-1}(8\delta
|m|+202\Delta_m)$ satisfies 
\[
\E\cro{H^2(\tilde{s},s)}\le\left[390\left(\inf_{m\in\M}\pa{H^2(s,S_{m})+
\pen(m)}+\frac{101\Sigma^2}{n}\right)+\varepsilon\right]\bigwedge2.
\]
\end{thm}
The only previous works on partition selection for histograms using squared Hellin\-ger loss
we know about are to be found in Castellan (1999 and 2000) and Birg\'e~(2006). Castellan's
approach is based on penalized maximum likelihood. This requires to make specific
restrictions on  the underlying density $s$, in particular that $s$ should be bounded away
from 0. For the problem of estimating a density on $\R$, her conditions on the family of
partitions are also more restrictive than ours since we can handle any countable families of
finite partitions into intervals. Nevertheless, in the multivariate case, our assumptions on the
partitions are more stringent. Birg\'e's approach based on aggregation of histograms built on
one half of the sample leads to more abstract but more general results. 

Let us now apply the above theorem to various families of models, systematically setting
$\pen(m)=n^{-1}(8\delta |m|+202\Delta_m)$ and $\varepsilon=n^{-1}$. We assume in this
section that $\lambda$ is the Lebesgue measure on $\X$.

\subsubsection{Example 1, continued}\label{U2a}
When $\X= [0,1)$, we use the family of models and weights of Section~\ref{A4a}. Our next
proposition shows that the PHE based on this simple family of models and weights has nice
properties for estimating various types of functions. The proof will be given in
Section~\ref{P1}.
%
\begin{prop}\label{P-ex3}
Let $\tilde{s}$ be the PHE based on the family of models and weights $\Delta_m$ defined in
Section~\ref{A4a}, $\varepsilon=n^{-1}$ and the penalty function
$\pen(m)=n^{-1}(8\delta |m|+202\Delta_m)$.

i) If $s\in{\S}_1(R)$, then
\begin{equation}
\E\left[H^2(\tilde{s},s)\right]\le  C\left\{\left[Rn^{-1}
\log\left(1+nR^2\right)\right]^{2/3}\vee n^{-1}\right\}.
\label{Ri2}
\end{equation}

ii) If $\sqrt{s}\in\S^2(w)$ where $w$ is a modulus of continuity  on $[0,1)$,  we define
$x_w$ to be the unique  solution of the equation $nxw^2(x)=1$ if $w(1)\ge n^{-1/2}$ and 
$x_w=1$ otherwise. Then
\begin{equation}
\E\left[H^2(\tilde{s},s)\right]\le  C(nx_w)^{-1}.
\label{Ri1}
\end{equation}
If, in particular, $\sqrt{s}$ belongs to the H\"older class $\H_\alpha^R$ with $R\ge
n^{-1/2}$, then $\E_s\left[H^2(\tilde{s},s)\right]\le
CR^{2/(2\alpha+1)}n^{-2\alpha/(2\alpha+1)}$.

iii) If $s\in\S^3(D,R)$ with $2\le D\le n$ and $R\ge2$, we get
\begin{equation}
\E_s\left[H^2(\tilde{s},s)\right]\le CDn^{-1}\log\left(nR/D\right).
\label{Ri3}
\end{equation}
\end{prop}
It is interesting to see to what extent the previous bounds (together with the trivial one,
$\E\left[H^2(\tilde{s},s)\right]\le2$, which always holds but which we did not include in
(\ref{Ri2}), (\ref{Ri1}) and (\ref{Ri3}) for simplicity) are optimal (up to the universal constants
$C$).  Many lower bounds on the minimax risk over various density classes are known for
classical loss functions. For squared Hellinger loss, some are given in Birg\'e (1983 and 1986)
and Birg\'e and Massart (1998). Many more are known for the squared $\IL_2$-loss, which
can easily be extended to squared Hellinger loss because their proofs are based on
perturbations arguments involving sets of densities for which both distances are equivalent. It
follows from these classical results that the bound we find for continuous densities are actually
optimal (see Birg\'e, 1983, p.211) while (\ref{Ri2}) is suboptimal because of the presence of
the $\log$ factor. We shall see below that the more sophisticated penalization strategy
introduced in Section~\ref{A4w} does solve the problem. The case of piecewise constant
functions is more complicated. If $D$ and the locations of the jumps were known, one could
use a single model corresponding to the relevant partition with $D$ intervals and get a risk
bound $CD/n$ corresponding to a parametric problem with $D$ parameters. Apart from the
constant $C$, this bound cannot be improved which shows that the study of uniform risk
bounds over $\S^3(D,R)$ is only of interest when $D\le n$ since otherwise a lower bound for
the risk is of the order of the trivial upper bound 2. When $D$ is smaller than $n$ the extra
$\log(nR/D)$ factor in (\ref{Ri3}) is due to the fact that we have to estimate the locations of
the jumps. The problem has been considered in Birg\'e and Massart (1998, Section~4.2 and
Proposition~2) where it is shown that a lower bound for the risk (when $n\ge 5D$ and
$D\ge9$) is $cDn^{-1}\log\left(nD^{-1}\right)$. Therefore our bound is optimal for moderate
values of $R$. We do not know whether the $\log R$ factor in the upper bound is necessary or
not.

\subsubsection{Improved risk bounds with a better weighting strategy}\label{U2b}
If we use the weights $\Delta_m$ defined in Section~\ref{A4w} to build $\tilde{s}$, we can
only improve (up to constants) the risk bounds given in Proposition~\ref{P-ex3} since the value
of $\Sigma$ does not change much while the new weights are not larger than the previous ones.
Besides, the values of  the weights have been substatially decreased for the partitions
belonging to $\M^1_T$. It turns out that piecewise constants functions on the elements of
$\M^1_T$ possess quite powerful approximation properties with respect to functions in Besov
spaces $B^\alpha_{p,\infty}([0,1))$ with $\alpha<1$ and monotone functions. These properties
are given in the following theorem which also includes the multidimensional case.
%
\begin{thm}\label{Besappro}
Let $\X=[0,1)^k$, $\M^k_T$ be the set of partitions $m$ of $\X$ defined in Section~\ref{A4c}
and, for $m\in\M^k_T$, let $S'_m$ be the cone $\left\{t=\sum_{I\in m}t_I\1_{I},
t_I\ge0\right\}$. For any $p>0$, $\alpha$ with $1>\alpha>k(1/p-1/2)_+$ and any function
$t$  belonging to the Besov space $B^\alpha_{p,\infty}([0,1)^k)$ with Besov semi-norm
$|t|_{B^\alpha_{p,\infty}}$, one can find some $t'\in\bigcup_{m\in\M^k_T}S'_m$ such
that
\begin{equation}
\|t-t'\|_2\le C(\alpha,k,p)|t|_{B^\alpha_{p,\infty}}|m|^{-\alpha/k},
\label{Eq-DeVo}
\end{equation}
where $\|\cdot\|_2$ denotes the $\Bbb{L}_2(dx)$-norm on $[0,1)^k$.

If  $t$ is a function of bounded variation on $[0,1)$, there exists
$t'\in\bigcup_{m\in\M^1_T}S'_m$  such that $\|t-t'\|_2\le C'\Var^*(t)|m|^{-1}$.
\end{thm}
The bound (\ref{Eq-DeVo}) is given in DeVore and Yu (1990). The proof for the bounded
variation case has been kindly communicated to the second author by Ron DeVore. With the
help of this theorem, we can now derive from Theorem~\ref{densite} the following improved
bounds the proof of which is straightforward.
%
\begin{prop}\label{P-ex6}
Let $\tilde{s}$ be the PHE based on the weights $\Delta_m$ defined in Section~\ref{A4w}. If
$\sqrt{s}$ is a function of bounded variation with $\Var^*\left(\sqrt{s}\right)\le R$ and in
particular if it belongs to $\S^1(R)$, then
\begin{equation}
\E\left[H^2(\tilde{s},s)\right]\le\min\left\{C(R/n)^{2/3},2\right\}\quad\mbox{for
}R\ge n^{-1/2}.
\label{Eq-BV}
\end{equation}
If  $\sqrt{s}\in B^\alpha_{p,\infty}([0,1))$ with $1>\alpha>(1/p-1/2)_+$ and
$\left|\sqrt{s}\right|_{B^\alpha_{p,\infty}}\le R$ with $R\ge n^{-1/2}$, then
\[
\E\left[H^2(\tilde{s},s)\right]\le\min\left\{C(\alpha,p)R^{2/(1+2\alpha)}
n^{-2\alpha/(1+2\alpha)},2\right\}.
\]
\end{prop}
It follows from classical lower bounds arguments that these bounds are minimax up to
constants. 

\subsubsection{The multidimensional case}\label{U2c}
When the density $s$ defined on $\X=\B_k$ can be written $s(x)=\Phi(\|x\|)$ for some
function $\Phi$ on $[0,1)$, we use the family of models introduced in Example~2. We then
obtain the risk bounds given in Propositions~\ref{P-ex3} and \ref{P-ex6} if we replace the
assumptions on $s$ by the same on $\Phi$. We omit the details.

If $\X=[0,1)^k$, $k\ge2$ and we use the family of models and weights described in
Section~\ref{A4c}, we get the following result.  
%
\begin{prop}\label{P-ex4}
Let $R\ge k^{-1}n^{-1/2}$. If $\sqrt{s}$ belong to $\H^R_{\alpha}([0,1)^k)$, then
\begin{equation}
\E\left[H^2(s,\tilde{s})\right]\le \min\left\{C(Rk)^{2k/(k+2\alpha)}
n^{-2\alpha/(2\alpha+k)},\,2\right\}.
\label{Eq-hold}
\end{equation}
More generally, if $\sqrt{s}$ belongs to $B^\alpha_{p,\infty}([0,1)^k)$ with
$1>\alpha>k(1/p-1/2)_+$ and $\left|\sqrt{s}\right|_{B^\alpha_{p,\infty}}\le R$, then
\[
\E\left[H^2(\tilde{s},s)\right]\le\min\left\{C(\alpha,k,p)R^{2k/(k+2\alpha)}
n^{-2\alpha/(k+2\alpha)},\,2\right\}.
\]
\end{prop}
%
\noindent{\em Proof:}
Let $m=\K_j$ be an element of $\M_R$. Then $\Delta_m=|m|=2^{kj}$ and the maximal
variation of a function of $\H^R_{\alpha}([0,1)^k)$ on an element of $m$ is bounded by 
$Rk2^{-j\alpha}$ so that $H^2(s,S_m)\le(Rk)^22^{-2j\alpha}$. It then follows from
Theorem~\ref{densite} that $\E\cro{H^2(\tilde{s},s)}\le C'
\left[(Rk)^22^{-2j\alpha}+n^{-1}2^{kj}\right]$. The lower bound on $R$ allows us to choose
$j\in\N$ such that $2^j\le\left(n(Rk)^2\right)^{1/(k+2\alpha)}< 2^{j+1}$ which leads to
\[
\E\cro{H^2(\tilde{s},s)}\le
C'\left[(Rk)^22^{2\alpha}\left(n(Rk)^2\right)^{-2\alpha/(k+2\alpha)}
+n^{-1}\left(n(Rk)^2\right)^{k/(k+2\alpha)}\right].
\]
The first bound follows since $2^{2\alpha}\le4$. The second bound can
be proved in the same way from (\ref{Eq-DeVo}).\cqfd
%
\subsection{Poisson processes}\label{U3}
Let us consider the stochastic framework corresponding to Problem~2 where $\nu$ is 
dominated by some given measure $M=\lambda$ on
$\X$ with density $s=d\nu/d\lambda$. This  implies that (\ref{Eq-fond}) holds as required. In this case, the
performances of the PHE $\tilde{s}$ are as follows. 
%
\begin{thm}\label{poisson}
Assume that the family $\M$ satisfies Assumption~{\bf H} and the weights $\ac{\Delta_m,\
m\in\M}$ are chosen so that \eref{som} holds. Then the estimator $\tilde{s}$ defined in
Section~\ref{A3} with $\pen(m)\ge3\delta|m|+6\Delta_m$ satisfies
\begin{equation}
\E\cro{H^2(\tilde{s},s)}\le390\left[\inf_{m\in\M}\cro{H^2(s,S_{m})+\pen(m)}+
3\Sigma^2\right]+\varepsilon.
\label{Eq-poiss}
\end{equation}
\end{thm}
This theorem should be compared with the results of Reynaud-Bouret~(2003) who uses more 
general families of projection estimators than just histograms based on partitions.
Nevertheless, for the problem we consider here, her choice of the $\IL^{2}$-loss induces some
restrictions on both the intensity and the collection of partitions at hand. For instance, the
intensity has to be bounded and the procedure requires some suitable estimation of its
sup-norm. As Castellan~(1999), she cannot deal with partitions with arbitrary small length.   

Let us now apply this theorem to our families of models, systematically setting
$\pen(m)=3\delta|m|+6\Delta_m$ and $\varepsilon=1$. In view of facilitating the
interpretation of the results to follow, it is convenient to use an analogy with density
estimation. This analogy, based on the following heuristics, allows to extrapolate the bounds
from one framework to the other.

We recall that observing the Poisson process $N$ of intensity $s$ is
equivalent to observing $\overline{N}$ i.i.d.\ random variables with density
$s'$, where $\overline{N}=N(\X)$ is a Poisson variable with parameter
$n=\int_{\X}s\,d\lambda$ and $s'=n^{-1}s$. With this in mind, and even though
$n$ need not be an integer, we can view the estimation of $s$ as an
analogue of the estimation of the density $s'$ from $n$ i.i.d.\ observations.
Pursuing into this direction, we may rewrite the risk in the Poisson case as
$\E\left[H^2(\tilde{s},s)\right]= n\E\left[H^2(n^{-1}\tilde{s},s')\right]$
and, setting $\tilde{s}_n=n^{-1}\tilde{s}$, view
$\E\left[H^2(\tilde{s}_n,s')\right]=n^{-1} \E\left[H^2(\tilde{s},s)\right]$ as
an analogue of the risk for estimating $s'$ from $n$ i.i.d.\ observations.
When $\sqrt{s}$ belongs to $\S^1(R)$, $\S^2(w)$ or $\S^3(D,R)$, then the
square-root of the density $s'=s/n$ belongs to $\S^1(Rn^{-1/2})$,
$\S^2(wn^{-1/2})$ or $\S^3(D,Rn^{-1/2})$ respectively (provided that $R^2\ge
n$ in the last case, since otherwise $\S^3(D,Rn^{-1/2})$ would not contain any
density). From these two remarks, we may conclude that a risk bound of the
form $f(R)$ in the Poisson case should be interpreted in the density case as
$n^{-1}f(Rn^{-1/2})$.

\subsubsection*{Example 1, continued}
Here we deal with a Poisson process $N$ on a finite interval of $\R$, which we may
assume, without loss of generality, to be $[0,1)$, of intensity $s$ with respect to the Lebesgue
measure $\nu$. To estimate $s$ we use the family of models of Example~1 with the
weights $\Delta_m$ defined in Section~\ref{A4w}. The resulting PHE $\tilde{s}$ has the
following properties which can be proved exactly like those given in
Propositions~\ref{P-ex3} and \ref{P-ex6}.
%
\begin{prop}\label{P-ex3p}
Let $w$ be a modulus of continuity  on $[0,1)$. We define $x_w$ to be the unique  solution of
the equation $xw^2(x)=1$ if $w(1)\ge1$ and  $x_w=1$ otherwise. Then
\begin{equation}
\E\left[H^2(\tilde{s},s)\right]\le  Cx_w^{-1}\quad\mbox{for all $s$ such that }
\sqrt{s}\in\S^2(w).
\label{Ri4}
\end{equation}
If, in particular, $\sqrt{s}$ belongs to the H\"older class $\H_\alpha^R$ with $R\ge1$,
then 
\[
\E\left[H^2(\tilde{s},s)\right]\le C R^{2/(2\alpha+1)}.
\]

Given $D\ge2$ and $R\ge2D$, we get
\begin{equation}
\E\left[H^2(\tilde{s},s)\right]\le CD\log(R/D)\quad\mbox{for all }s\in\S^3(D,R).
\label{Ri5}
\end{equation}

If $\sqrt{s}$ belongs to $\S^1(R)$ with $R\ge 1$, then $\E\left[H^2(\tilde{s},s)\right]\le
CR^{2/3}$.

If  $\sqrt{s}\in B^\alpha_{p,\infty}([0,1))$ with $1>\alpha>(1/p-1/2)_+$ and
$\left|\sqrt{s}\right|_{B^\alpha_{p,\infty}}\le R$ with $R\ge 1$, then
$\E\left[H^2(\tilde{s},s)\right]\le CR^{2/(1+2\alpha)}$.
\end{prop}
For the sake of simplicity, let us assume that $n=\int_{\X} sd\lambda$ is an integer. The
connection established above between the estimation of a density and that of the intensity of a
Poisson process shows that Proposition~\ref{P-ex3p} is actually a perfect analogue of
Propositions~\ref{P-ex3} and~\ref{P-ex6}.  Namely, when $\sqrt{s}$ belongs to $\S^1(R)$ or
$\S^2(w)$ or $s\in\S^3(D,R)$ and $s'=s/n$ then $\sqrt{s'}$ respectively belongs to
$\S^1(Rn^{-1/2})$ or $\S^2(wn^{-1/2})$ or $s'\in\S^3(D,Rn^{-1})$ and the risk bounds we get
for estimating the intensity $s$ (with respect to the $H^2/n$-loss) are the same as those obtained
from a $n$ sample for estimating the density $s'$ (with the $H^2$-loss).

\subsubsection*{Example 2, continued}
If we observe a Poisson process on $\X=\B_k$ with intensity  $s(x)=\Phi(\|x\|)$ with respect
to the Lebesgue measure for $\Phi$ some function on $[0,1)$ and consider the family of
models introduced in Example~1 we obtain the risk bounds given in Proposition~\ref{P-ex3p}
if we replace the assumptions on $s$ by the same on $\Phi$.

\subsubsection*{Example 3, continued}
If $\X=[0,1)^k$ with $k\ge2$, we use the models and weights defined in Section~\ref{A4c}. 
Proceeding as for Proposition~\ref{P-ex4} we get:
%
\begin{prop}\label{P-ex5}
Let $\sqrt{s}$ belong to $\H^R_{\alpha}([0,1)^k)$, then
\[
\E_s\left[H^2(s,\tilde{s})\right]\le C(Rk\vee1)^{2k/(k+2\alpha)}.
\]
If $\sqrt{s}$ belongs to $B^\alpha_{p,\infty}([0,1)^k)$ with $1>\alpha>k(1/p-1/2)_+$ and
$\left|\sqrt{s}\right|_{B^\alpha_{p,\infty}}\le R$, then
\[
\E\left[H^2(\tilde{s},s)\right]\le C(\alpha,k,p)(R\vee1)^{2k/(k+2\alpha)}.
\]
\end{prop}
As shown by the proof of Proposition~\ref{P-ex4}, we only use the partitions in $\M_R$ to get
(\ref{Eq-hold}) so that it would be of little use to introduce other partitions if we only wanted to
estimate intensities such that $\sqrt{s}$ belong to $\H^R_{\alpha}([0,1)^k)$. The interest of
considering the larger family $\M$ and to have a special definition of $\Delta_m$ when
$m\in\M_T$ is that it allows to improve the results when we deal with less regular functions
than those for which $\sqrt{s}$ belong to $\H^R_{\alpha}([0,1)^k)$, in particular those
functions that belong to Besov spaces $B^\alpha_{p,\infty}([0,1)^k)$ with $1>\alpha>k/p$.
To illustrate this fact, let us study the estimation of those intensities $s$ such that $\sqrt{s}$
has the following specific structure. Given the nonempty set $\V$ which is a finite union of
elements of  $\K$, there is a smallest integer $\bar{\j}$ such that $\V$ can be written as the
union of $N$ elements of $\K_{\bar{j}}$ with a volume $V=N2^{-k\bar{\j}}>0$. To avoid
trivialities, we assume that $\bar{\j}>0$, hence $V<1$. 
%
\begin{prop}\label{P-spepar}
Let $s$ be an intensity on $[0,1)^k$ such that $\sqrt{s}\1_{\V}$ belongs to
$\H^R_{\alpha}(\V)$ with $R\ge1$ while $\sqrt{s}\1_{\V^c}$ is constant and let $\tilde{s}$
be the PHE based on the weights $\Delta_m$ defined in Section~\ref{A4c}.  Then
\begin{equation}
\E\cro{H^2(\tilde{s},s)}\le C\inf_{m\in\M}B_m\quad\mbox{with }B_m=
H^2(s,S_m)+|m|+\Delta_m
\label{Eq-Bm}
\end{equation}
and
\begin{eqnarray}
B_m&\le&C\min\left\{2^{k\bar{\j}}+V^{k/(k+2\alpha)}(kR)^{2k/(k+2\alpha)}\right. ;
\label{mi1}\\&&\mbox{}\hspace{13mm}V\left[k\bar{\j}2^{k\bar{\j}}+(kR)^{2k/(2\alpha+k)}
\left[\log\left(Rk\right)\right]^{2\alpha/(2\alpha+k)}\right];\label{mi2}\\&&\mbox{}
\hspace{12mm}\left.V\left[2^{k}\bar{\j}2^{k\bar{\j}}+(kR)^{2k/(2\alpha+k)}\right]\right\}. 
\label{mi3}
\end{eqnarray}
\end{prop}
%
\noindent{\em Proof:}
Since (\ref{Eq-Bm}) is merely a consequence of Theorem~\ref{poisson} with the choice
$\pen(m)=3\delta|m|+6\Delta_m$ and $\varepsilon=1$, we only have to bound $B_m$. Let
us first consider a regular partition $m=\K_j$. If $j<\bar{\j}$, the bias $H^2(s,S_m)$ may be
arbitrarily large since the intensity $s$ may be arbitrarily large on $\V$ while it may be small on
$\V^c$. For $j\ge\bar{\j}$, the argument used for the proof of Proposition~\ref{P-ex4} shows
that on $\V$, $\sqrt{s}$ can be approximated uniformly by an element of $S_m$ with a precision
at least $Rk2^{-j\alpha}$ so that $H^2(s,S_m)\le VR^2k^22^{-2j\alpha}$ and $B_m\le
VR^2k^22^{-2j\alpha}+2^{kj+1}$. If  $\left[VR^2k^2\right]^{1/(2\alpha+k)}\le2^{\bar{\j}}$ we
set $j=\bar{\j}$ and otherwise choose $j$ so that $2^j\le\left[VR^2k^2\right]^{1/(2\alpha+k)}
<2^{j+1}$. This leads to (\ref{mi1}).
 
If we set $m=m_p\vee\K_0$ with $p$ being the set of those $N2^{k(j-\bar{\j})}=V2^{kj}\ge1$
elements of  $\K_j$ ($j\ge\bar{\j}\ge1$) that exactly cover $\V$, we get, since $k\ge2$
\[ 
B_m\le VR^2k^22^{-2j\alpha}+(kj+1)V2^{kj}+1\le Vk\left[R^2k2^{-2j\alpha}+2j2^{kj}\right]. 
\]
If  $\left[k^2R^2/\log(kR)\right]^{1/(2\alpha+k)}<2^{\bar{\j}}$ we set $j=\bar{\j}$ and
otherwise choose $j$ so that $2^j\le\left[k^2R^2/\log(kR)\right]^{1/(2\alpha+k)}<2^{j+1}$
which finally leads to (\ref{mi2}).

To study the approximation properties of the elements of $\M^k_T$ let us consider a particular
cube $K'=K_{\bar{\j},\bm{l}}\in\V\cap\K_{\bar{\j}}$. Identifying the partitions in $\M^k_T$
with the trees from which they derive, we can design an element $m_{K'}$ of $\M^k_T$ with
$2^k-1$ terminal nodes at each level 1 to $\bar{\j}$ and the remaining node $K'$ at level
$\bar{\j}$. Then we keep only non-terminal nodes up to level $j\ge\bar{\j}$, all nodes at this last
level $j$ being terminal, so that their number is  $2^{k(j-\bar{\j})}$. The total number of terminal
nodes of the tree is therefore $\bar{\j}(2^k-1)+2^{k(j-\bar{\j})}$. We can repeat this operation for
each of the $N$ cubes in $\V\cap\K_{\bar{\j}}$ keeping the value of $j$ fixed. This results in $N$
similar trees. We finally consider the smallest complete tree $m$ that contains the $N$ previous
ones. Its number of terminal nodes is then bounded by
$N\left[\bar{\j}(2^k-1)+2^{k(j-\bar{\j})}\right]$ so that
\[
B_m\le V(Rk)^22^{-2j\alpha}+2N\left[\bar{\j}(2^k-1)+2^{k(j-\bar{\j})}\right]
\le2V\left[R^2k^22^{-2j\alpha}+\bar{\j}2^{k(\bar{\j}+1)}+2^{kj}\right]. 
\] 
If  $\left(k^2R^2\right)^{1/(2\alpha+k)}<2^{\bar{\j}}$ we set $j=\bar{\j}$ and
otherwise choose $j$ so that $2^j\le\left(k^2R^2\right)^{1/(2\alpha+k)} <2^{j+1}$,
which leads to (\ref{mi3}).\cqfd\vspace{2mm}\\
A comparison of the three bounds (\ref{mi1}), (\ref{mi2}) and (\ref{mi3}) shows that
(\ref{mi3}) is always better if we omit the influence of $\bar{j}$ and $k$ but the situation
becomes more involved if we take into account the effect of $k$ and $\bar{j}$. Depending on the
values of $V,R,\bar{\j},\alpha$ and $k$, each type of partition may be the best which justifies
to introduce them all. 

\noindent{\em Remark:} An analogue of Proposition~\ref{P-spepar} holds for density
estimation.

%
\subsection{Non-negative random vectors}\label{U1}
Let us recall from the introduction that we observe an $n$-dimensional random vector with
independent nonnegative components $N_1,\ldots,N_n$ and respective distributions
depending on positive parameters $s_1,\ldots,s_n$. One should think of the $N_i$ as Poisson or
binomial random variables with unknown expectations $s_i$. More generally, we assume that
there exist some known constants $\kappa>0$ and $\tau\geq 0$ such that for all
$i\in\X=\ac{1,\ldots,n}$
\begin{equation}
\log\pa{\E\cro{e^{z\pa{N_i-s_i}}}}\le \kappa\frac{z^2s_i}{2(1-z\tau)}\quad
\mbox{for all }z\in\left[0,\frac{1}{\tau}\right[,
\label{eq:lap1e}
\end{equation}
with the
convention $1/\tau=+\infty$ if  $\tau=0$, and
\begin{equation}
\log\pa{\E\cro{e^{-z\pa{N_i-s_i}}}}\le\kappa \frac{z^2s_i}{2} \ \quad\mbox{for
all }z\geq 0.
\label{eq:lap2e}
\end{equation}
In the case of Poisson or binomial random variables, one can take $\kappa=\tau=1$ as we
shall see below. 

Our aim is to estimate the function $s$ from $\X$ to $\R_+$ given by $s(i)=s_i$. Here we
denote by $\lambda$ the counting measure on $\X$ and set $Y\equiv1$. Hence $M=\lambda$
and
$N(A)=\sum_{i\in A}N_i$. Then $\LL$ can be identified with $\R_+^n$,
$\E\cro{N(A)}=\int_{A} sd\lambda$ as required and $H^2(t,t')=\sum_{i=1}^n
\cro{\sqrt{t(i)}-\sqrt{t'(i)}}^2$ for $t,t'\in\LL$.
\begin{thm}\label{cas-discret}
Assume that (\ref{eq:lap1e}) and (\ref{eq:lap2e}) hold, that the family $\M$ satisfies
Assumption~{\bf H} and the weights $\ac{\Delta_m,\ m\in\M}$ are chosen so that
\eref{som} holds. Let
$\pen(m)\ge\kappa\left[\delta\pa{1+K^2}|m|+3K^2\Delta_m\right]$ with
\[ 
K=\sqrt{2}\quad\mbox{if }\tau\le\kappa;\qquad K=\frac{\sqrt{2}}{2} 
+\sqrt{\frac{\tau}{\kappa}-\frac{1}{2}}\quad\mbox{if }\tau>\kappa;
\] 
and let $\tilde{s}$ be the PHE defined in Section~\ref{A3}. Then 
\[
\E\cro{H^2(\tilde{s},s)}\leq390\left[\inf_{m\in\M}
\cro{H^2(s,S_{m})+\pen(m)}+(3/2)\kappa K^2\Sigma^2\right]+\varepsilon.
\]
\end{thm}
Let us first check that some classical distributions do satisfy Inequalities~\eref{eq:lap1e}
and~\eref{eq:lap2e}. If $N_i$ is a binomial random variable with parameters $n_i,p_i$ then for
all $z\in\R$, 
\begin{equation}
\log\left(\E\cro{e^{z(N_i-s_i)}}\right)\le s_i\pa{e^{z}-z-1}\quad\mbox{with}\quad
s_i=n_ip_i.
\label{Eq-bin}
\end{equation}
If  $N_i$ is a Poisson random variable with parameter $s_i$, then equality holds in
(\ref{Eq-bin}). Using the bounds $e^{z}-z-1\le z^2/[2(1-z)]$ for $z\in[0,1[$ and
$e^{z}-z-1\le z^2/2$ for $z<0$ we derive that, in both cases, \eref{eq:lap1e}
and~\eref{eq:lap2e} hold with $\kappa=\tau=1$. If $N_i$ has a Gamma distribution
$\Gamma(s_i,1)$, $\E\cro{N_i}=s_i$ and, following the proof of Lemma~1 of Laurent and
Massart (2000), we deduce that ~\eref{eq:lap1e} and~\eref{eq:lap2e} hold again with
$\kappa=\tau=1$.  More generally, it follows from some version of Bernstein's Inequality --- see 
Lemma~8 of Birg\'e and Massart (1998) --- that~\eref{eq:lap1e} holds as soon as 
\[
\E\left[(N_i)^p\right]\le\kappa\frac{p!}{2}s_i \tau^{p-2},
\qquad\mbox{for all }i\in\X\quad\mbox{and}\quad p\geq 2.
\]
Inequality~\eref{eq:lap2e} is always satisfied if  $N_i\leq\kappa$. Indeed it follows from
$$e^{-zx}\leq 1-zx+z^2x^2/2,\ \  \forall x,z\geq 0$$
that all non-negative random variables $X$
bounded by $\kappa$ satisfy
\[
\E\cro{e^{-zX}}\leq 1-z\E[X]+\frac{z^2\E[X^2]}{2}\le\exp\pa{-z\E[X]+\kappa
  z^2\E[X]/2}.
\]
The results of Kolaczyk and Nowak (2004), which are based on some sort of discretized
penalized maximum likelihood estimator in the spirit of Barron and Cover (1991), have some
similarity with ours but they assume that the components of the vector $s$ belong to some
known interval $[c,C]$, $c>0$ and they explicitely use the values of $c$ and $C$ in the
construction of their estimator. Such an assumption, which implies, as in the case of density
estimation, that squared Hellinger distance and Kullback divergence are equivalent also greatly
simplifies the estimation problem.

\subsubsection*{Example 4, continued}
Setting
\begin{equation}
\pen(m)=\kappa\left[\pa{1+K^2}|m|+3K^2\Delta_m\right]\qquad\mbox{and}\qquad
\varepsilon=1.
\label{Eq-ex1}
\end{equation}
and using $\log\binom{n-1}{D-1}\le(D-1)(1+\log[(n-1)/(D-1)])$ with the convention
$0\log((n-1)/0)=0$ we get the risk bound
\begin{equation}
\E\cro{H^2(\tilde{s},s)}\le C(\kappa,K)\inf_{m\in\M}\left\{H^2(s,S_m)+|m|
+(|m|-1)\log\left(\frac{n-1}{|m|-1}\right)\right\}.
\label{eq-78}
\end{equation}
If, for instance, $s$ itself belongs to some $S_m$ with a small value of $|m|$, which
corresponds to a piecewise stationary process $(N_i)_{1\le i\le n}$ with a few distribution
changes, the risk is bounded by $C(\kappa,K)|m|\log n$.

Another interesting situation corresponds to the case of a monotone sequence $(s_i)_{1\le
i\le n}$, i.e.\ a monotone function $s$ on $\X$ that we may assume, without loss of
generality to be nondecreasing.
%
\begin{prop}\label{P-mon6}
Let the sequence $s_i, 1\le i\le n$ be nondecreasing with $\sqrt{s_n}-\sqrt{s_1}=R$, then 
the PHE $\tilde{s}$ based on the models of Example~4 with $\pen$ and $\varepsilon$
given by (\ref{Eq-ex1}) satisfies the following risk bounds with a constant $C$ depending
only on $\kappa$ and $K$:
\begin{itemize}
  \item if  $R^2\le n^{-1}\log n$, then $\E\cro{H^2(\tilde{s},s)}\le
C(\kappa,K)\left(nR^2+1\right)$;
  \item if $R\ge n/\sqrt{3}$, then $\E\cro{H^2(\tilde{s},s)}\le C(\kappa,K)n$;
                \item otherwise
$\E\cro{H^2(\tilde{s},s)}\le C(\kappa,K)\left[R\sqrt{n}\log(n/R)\right]^{2/3}$.
\end{itemize}
\end{prop}
%

\noindent{\em Remark:} If we restrict ourselves to the case $n=2^k$, we can turn any 
function $s$ on $\X$ into a function $s'$ on $[0,1)$ by setting
$s'=\sum_{i=1}^ns(i)\1_{[(i-1)2^{-k},i2^{-k})}$. This transformation will, in particular, 
preserve the monotonicity properties of the functions. One could then estimate $s'$ using
the more sophisticated families of weights that we introduced in Section~\ref{A4w}. The use of
this strategy would improve the estimation of monotone functions, removing the logarithmic
factors. 

\subsubsection*{Example 5, continued} 
 Choosing $\pen$ and $\varepsilon$ as in (\ref{Eq-ex1}) and using the same
arguments as for Example~1, we derive an analogue of (\ref{eq-78}) with $n$ replacing $n-1$
in the logarithmic factor. If we assume  that $s_i=\overline{s}$ for $i\not\in I$ with $|I|=k$,
then
$H^2(s,S_m)=0$ for some
$m\in\M_k$ and 
\[
\E\cro{H^2(\tilde{s},s)}\le C(\kappa,K)[k+1+k\log(n/k)].
\]

%
%
\section{Special counting processes on the line}\label{F}
Let $\X$ be some interval of $\R_+$ of the form $[0,\zeta)$ where $0<\zeta\le+\infty$ with
its Borel $\sigma$-algebra $\A$.  We recall that a (univariate) counting process $\widetilde N$
on
$\X$ is a cadlag (right-hand continuous and left-hand limited) process from $\X$ to $\R_+$,
vanishing at time $t=0$, with piecewise constant and nondecreasing paths having jumps of size
$+1$ only.  The use of counting processes in statistical modeling is developed in great details in
the book by Andersen {\it et al.}~ (1993) where the interested reader will find many concrete
situations for which these processes naturally arise. Typically, $\widetilde{N}_t$ counts the
number of occurrences of a certain event from time 0 up to time $t$. The {\it jumping times} of the
process give the dates of occurrence of the event. A counting process can be associated to a
random measure $N$ on $\X$ whose cumulative distribution function is the counting
process itself, i.e.\ $N([0,t])=\widetilde N_t$ for all $t\in \X$. In the sequel, we shall not
distinguish between the counting process $\widetilde N$  and its associated measure $N$. 

In this paper, we consider a phenomenon which is described by some bounded counting
process $N^*$ on $\X$ such that $N^*(\X)\leq k$ a.s.\ for some known integer $k$. This
means that $N^*$ describes an event that occurs at most $k$ times during the period $\X$. We
also assume that there exist  a deterministic measure $\lambda$ on $\X$, a deterministic
nonnegative function $s\in\IL_1(\X,d\lambda)$ and  a nonnegative observable process
$Y^*$ bounded by 1 on $\X$ such that
\begin{equation}\label{aalen}
  \E\cro{N^*([0,t])}=\E\cro{\int_0^t sY^* d\lambda}\quad\mbox{for all }t\in\X.
\end{equation}
We actually observe an aggregated counting process $N$ which is the sum of
$n$  i.i.d.\ processes $N^j,\ j=1,\ldots,n$ with the same distribution as $N^*$. The fact that the
measure $N^j$ is determined by its cumulative distribution function and (\ref{aalen}) imply
that there are i.i.d.\ observable processes $Y^j$, $j\in\ac{1,\ldots,n}$ with the distribution of
$Y^*$ such that 
\[
\E\cro{N^j(A)}=\E\cro{\int_AsY^j d\lambda}\quad\mbox{for all }A\in\A
\quad\mbox{and}\quad 1\le j\le n.
\]
Therefore (\ref{Eq-fond}) holds with $M=Yd\lambda$ and $Y=\sum_{j=1}^n Y^j$. For such
counting processes, we can prove the following result.
\begin{thm}\label{cp}
Assume that there exist a positive integer $k$ and a positive number $\kappa'$, both known,
such that $N^*(\X)\le k$ a.s., (\ref{aalen}) holds and ${\rm Var}\cro{\int_{I}sY^*d\lambda}
\le\kappa'\E\cro{\int_I sY^*d\lambda}$ for all intervals $I\subset\X$. Assume moreover
that $\int_{\X} sd\lambda< +\infty$ and the aggregated process $N$ satisfies
(\ref{Eq-fond1}). Let us choose a family $\M$ satisfying Assumption~{\bf H} and weights
$\ac{\Delta'_m,\ m\in\M}$ such that 
\begin{equation}
\sum_{m\in\M}\exp[-\eta\Delta'_m]=\Sigma'(\eta)<+\infty\quad\mbox{for }
\eta=k\left(k+\int_{\X} s d\lambda\right)^{-1}.
\label{Eq-A3}
\end{equation}
Then the estimator $\hat s_{\hat m}$ defined in Section~\ref{A3} with
$\pen(m)\ge16\delta|m|(k+\kappa')+404k\Delta'_m$ satisfies
\begin{eqnarray*}
\lefteqn{\E\cro{H^2(\hat{s}_{\hat m},s)}}\hspace{15mm}\\&\le&
390\left(\E\cro{\inf_{m\in\M}\pa{H^2(s,S_{m})+\pen(m)}}+
404k\eta^{-1}[\Sigma'(\eta)]^2\right)+\varepsilon\\&\le&
390\left(\inf_{m\in\M}\left\{\E\cro{H^2(s,S_{m})}+\pen(m)\right\}+
404k\eta^{-1}[\Sigma'(\eta)]^2\right)+\varepsilon.
\end{eqnarray*}
\end{thm}
In the last bound, $\E\cro{H^2(s,S_{m})}$ plays the role of a bias term which can be bounded in
the following way. Let us set
\[ 
S'_m=\left\{t=\sum_{I\in m\cap \J}t_I\1_{I}\quad\mbox{with }
t_I\ge0\mbox{ for all }I\in m\right\}\bigcap\LL,
\] 
where the $t_I$ are now deterministic. Then $S'_m\subset S_m$, hence $H^2(s,S_m)\le H^2(s,
S'_m)$ and, for $t\in S'_m$,
\[
H^2(s,t)=\int_{\X}\left(\sqrt{s}-\sqrt{t}\right)^2Y\,d\lambda\le n
\int_{\X}\left(\sqrt{s}-\sqrt{t}\right)^2\,d\lambda,
\]
since $Y\le n$. Finally 
\[
\Bbb{E}\left[H^2(s,S_m)\right]\le n\inf_{t\in S'_m}
\int_{\X}\left(\sqrt{s}-\sqrt{t}\right)^2\,d\lambda=b_{m}^{2}(s)
\]
and
\[
\E\cro{H^2(\hat{s}_{\hat m},s)}\le 390\left(\inf_{m\in\M}\left\{b_{m}^{2}(s)+\pen(m)\right\}+{404k\over \eta}[\Sigma'(\eta)]^2\right)+\varepsilon.
\]
Note that the present framework includes, as a particular case, density estimation, if we observe
an $n$-sample $X_1,\ldots,X_n$ with density $s$ with respect to $\lambda$ and set
$N^j(A)=\1_{A}(X_j)$. Then $Y=n$ and $H^2(s,t)=n\int_{\X}\left(\sqrt{s}-\sqrt{t}\right)^2
\,d\lambda$ which corresponds to using the distance $H$ of Section~\ref{U2} multiplied by
$\sqrt{n}$. Up to this scaling factor, the previous risk bound is analogue to that for estimating
densities we get in Theorem~\ref{densite}.

In order to derive risk bounds which are similar to those given in Proposition~\ref{P-ex3}, we
have to distinguish between two situations. The most favorable one occurs when we know an
upper bound $\Gamma$ for $\int_{\X}s\,d\lambda$, in which case, since  $0\le Y^*\le1$,
\[
{\rm Var}\cro{\int_{I}sY^*d\lambda}\le\E\cro{\left(\int_I sY^*d\lambda\right)^2}
\le\left(\int_{\X}s\,d\lambda\right)\E\cro{\int_I sY^*d\lambda}
\]
and we can set $\kappa'=\Gamma$. Moreover, assuming that (\ref{som}) holds,
we can choose $\Delta'_m=\left(1+k^{-1}\Gamma\right)\Delta_m$ without any
further restriction on the family of models. Using the same family of partitions as in the density
case, we recover the bounds of Propositions~\ref{P-ex3} and \ref{P-ex6} up to the factor $n$
corresponding to the rescaling of the distance $H$.

Let us now turn to the less favorable situation where no bound for $\int_{\X}s\,d\lambda$ is known, which is the typical case for Problem~4. As we shall see the number $\kappa'$ can still be computed. As to (\ref{Eq-A3}) it will be satisfied with
$\Delta'_m=|m|$ as soon as the number of models such that $|m|=D$ is bounded
independently of $D$. Restricting ourselves to the family  $\M_R$ of regular partitions, we
recover, up to the factor $n$, the bounds provided by case ii) of Proposition~\ref{P-ex3}. 

\subsection{Survival analysis with right-censored data\label{F1}}
Let us now consider the framework of Problem~3, denoting by $P_{T}$ the common 
distribution of the $T_{i}$.  We consider the counting process $N$ on $\R_+$ defined by
$N=\sum_{j=1}^n N^j$ where $N^j(A)= \1_{\{\widetilde{T}_j\in A,\,D_j=1\}}$ for all
measurable subsets $A$ of $\R_+$, so that we can take $k=1$. Then the variables $N^j(A)$,
$1\le j\le n$ are i.i.d.\ Bernoulli random variables. We define $s$ to be the hazard rate of the
survival times, i.e.\ $s(t)=p(t)/\P[T_1\ge t]$ for $t>0$. Since $s$ is not integrable on $\R_+$
we shall restrict ourselves to some bounded interval $\X$ of $\R_+$, which we can take,
without loss of generality, to be $[0,1)$ if we assume that $\P[T_1\geq 1]>0$. We also assume
here that the censorship satisfies for all $t\ge 0$,
\begin{equation}\label{eq-compt}
\E\cro{N^j([0,t])}=\E\cro{\int_{0}^t s(u)Y^j(u)du},\quad\mbox{with }
Y^j(t)=\1_{\widetilde{T}_{j}\ge t},
\end{equation}
which means that (\ref{aalen}) holds. Equality~\eref{eq-compt} is clearly satisfied when
$C_j=T_j$ for all $j$, i.e.\ when the data are uncensored. It is also satisfied when the censorship
is independent of the survival time, i.e.\ when $C_j$ and $T_j$ are independent for all $j$.
Indeed, we then have for all $j$ and $t\ge0$, by Fubini Theorem and independence,
\begin{eqnarray*}
\E\cro{\int_{0}^t s(u)Y^j(u)du}&=&\E\cro{\int_{0}^t \frac{p(u)}{\P(T_j\geq u)}
\1_{C_j\geq u}\1_{T_j\geq u}du}\\&=&\int_{0}^t \frac{p(u)\P(T_j\geq u)\P(C_j\geq u)}
{\P(T_j\geq u)}du\\&=& \int\1_{[0,t]}(u) \P(C_j\geq u)dP_T(u)\\
&=& \P\cro{T_j\leq t, T_j\leq C_j}= \E\cro{N^j([0,t])}.
\end{eqnarray*}
%
\begin{prop}\label{P-surv}
If the processes $N^j$ satisfy (\ref{eq-compt}), the assumptions of Theorem~\ref{cp} hold
with $k=1$, $\kappa'=2$ and $\int_{\X} sd\lambda=-\log(\P[T_1\ge1])$. 
\end{prop}
From a practical point of view, one can always estimate $\P[T_1\geq 1]$ accurately enough to
assume that an upper bound $\Gamma$ for $\int_{\X} sd\lambda$ is known. We can
therefore apply Theorem~\ref{cp} to the the family of models of Example~1 with the weights
$\Delta_m$ given in Section~\ref{U2}, setting $\Delta'_m=(1+\Gamma)\Delta_m$. We then
obtain perfect analogues of Propositions~\ref{P-ex3} and~\ref{P-ex6} with constants $C$ now 
depending on $\Gamma$. To avoid redundancy, we leave the precise statement of the risk
bounds to the reader. 

\subsection{Transition intensities of Markov processes}\label{F2}
Within the framework of Problem~4, we associate to $T_{1,0}$ the counting process $N^*$ 
defined for $t\geq 0$ by
$N^*([0,t])=\1_{\{T_{1,0}\le t\}}$ so that 
\begin{equation}
\E\cro{N^*([0,t])}=\int_0^t p(u)du=\E\cro{\int_{0}^t \1_{\{X_{u-}=1\}}s(u)du}
\label{eq:mark}
\end{equation}
and (\ref{aalen}) holds with $Y^*(u)=\1_{\{X_{u-}=1\}}$.
Our aim here is to estimate $s$ on some bounded interval $\X$ of $\R_+$ from the
observation of the counting process $N=\sum_{j=1}^n N^j$ where the $N^j$'s are i.i.d.\
copies of $N^*$ associated to $n$ i.i.d.\ copies $X^1,\ldots,X^n$ of the process $X$. If $X$
takes only the two values 0 and 1 and a.s.\ starts from 1 to reach 0, then the problem reduces to
estimating the density $p$ of $T_{1,0}$; it becomes novel when we have at least three states. In
any case, we  get the following result.
%
\begin{prop}\label{P-mark}
If  the weights $\Delta'_m$ satisfy $\sum_{m\in\M}
\exp[-\eta\Delta'_m]<+\infty $ for all $\eta>0$ and $\int_{\X}s(t)dt<+\infty$ then Theorem~\ref{cp} applies with $k=1$ and
$\kappa'=2$. 
\end{prop}
\section{A unifying result}\label{B}
We want here to analyze our estimation procedure from the general point of view described in
Section~\ref{A} and prove a risk bound for the estimator $\tilde{s}$, from which we shall
be able to derive the previous risk bounds corresponding to all the specific frameworks
that we considered. For this we introduce the following approximation for $s$ in $S_m$:
\begin{equation}
\overline{s}_m=\sum_{I\in m\cap \J}\frac{s_I}{\lambda(I)}\1_{I}
\quad\mbox{with }s_I=\int_Is\,d\lambda.
\label{eq:def-bars}  
\end{equation}
We need here a bound for $H^2\left(\hat{s}_m,\overline{s}_m\right)$ which holds
uniformly for $m\in\overline{\M}$. It takes the following form:\vspace{2mm}\\
{\bf H'}\,:\, There exist three positive constants $a,b$ and $c$, $c\ge1$ such that,
for any $m\in\overline{\M}$,
\begin{equation}
\P\cro{H^2(\hat s_m,\overline{s}_m)\ge  c|m|+bz}\le a\exp[-z]
\quad\mbox{for all }z\ge0\vspace{2mm}.
\label{debase}
\end{equation}
We can now derive bounds for the risk of  the estimator $\tilde{s}$ defined in
Section~\ref{A3}.
\begin{thm}\label{main}
Let Assumptions~{\bf H} and {\bf H'} hold and the weights $\Delta_m$ satisfy
(\ref{som}). Let the penalty $\pen(m)$ be given by 
\begin{equation}
\pen(m)\ge c\delta|m|+b\Delta_m.
\label{Eq-pen}
\end{equation}
and $\hat{m}$ be any element of  $\M$ satisfying (\ref{Eq-def1}). Then the 
estimator $\tilde{s}=\hat{s}_{\hat m}$ satisfies 
\begin{equation}
\E\cro{H^2(\tilde{s},s)}\le390\left(\E\cro{\inf_{m\in\M}
\pa{H^2(s,S_{m})+\pen(m)}}+ ab\Sigma^2/2\right)+\varepsilon.
\label{Eq-risk1}
\end{equation}
\end{thm}
Note that such a result has been obtained without any assumption on the underlying space
$\X$ and the true value $s$ of the parameter, apart from the fact that it belongs to $\LL$. Note also
that in (\ref{Eq-risk1}), the infimum over $m\in\M$ occurs inside the expectation, which makes a
difference when $M$, and therefore $H(s,S_m)$, is random.

As we have previously seen, $\delta\le2$ for all the models we consider. Moreover, we shall see
in Sections~\ref{P1a}, \ref{P3a} and \ref{P4t} that for Problems~0, 1 and 2, $a=1$ and $b$ and
$c$ take the form $b=b'C_P$ and $c=c'C_P$ where $b'$ and $c'$ are numerical constants and
$C_P$ depends of the problem we consider (for instance $C_P=n^{-1}$ for density estimation). If
we  choose $\pen(m)=c_0C_P(|m|+\Delta_m)$ for some suitable numerical constant $c_0$ and
$\varepsilon\le C_P$, it follows that (\ref{Eq-risk1}) becomes
\begin{eqnarray*} 
\lefteqn{\E\cro{H^2(\tilde{s},s)}}\hspace{10mm}&&\\
 &\le& 390\left(\E\cro{\inf_{m\in\M}
\pa{H^2(s,S_{m})+c_0C_P(|m|+\Delta_m)}}+ 2b'C_P\Sigma^2/2\right)+C_P,
 \end{eqnarray*}
which gives (\ref{Eq-risk8}). If there is only one model $m$ in the family $\M$, we can fix
$\Delta_m=0$, hence $\Sigma=1$, which leads to (\ref{Eq-risk9}).
\begin{proof}
Let $m^*$ be an arbitrary element of $\M$.  It follows from the definition of $\D$ that for
any
$m\in\M$, $H^2(\hat{s}_m,\hat{s}_{m^*})\le\D(m)\vee\D(m^*)$. Therefore,
\begin{equation}
H^2(\hat{s}_{\hat{m}},\hat{s}_{m^*})\le\D(\hat{m})\vee\D(m^*)\le
\D(m^*)+\varepsilon/3,
\label{Eq-A6}
\end{equation}
by (\ref{Eq-def1}). It also follows from (\ref{Eq-def2}) that, if $T_{m,m^*}\le0$, then
\begin{equation}
H^2(\hat{s}_{m},\hat{s}_{m\vee m^*})-H^2(\hat{s}_{m^*},\hat{s}_{m\vee m^*})\le
16[\pen(m^*)-\pen(m)].
\label{Eq-A1}
\end{equation}
Moreover
\begin{eqnarray*}
\lefteqn{H^2(\hat{s}_{m},\hat{s}_{m\vee m^*})-H^2(\hat{s}_{m^*},\hat{s}_{m\vee m^*})}
\hspace{20mm}\\&=&\int\hat{s}_m\,d\lambda-\int\hat{s}_{m^*}\,d\lambda
+2\int\pa{\sqrt{\hat{s}_{m^*}}-\sqrt{\hat{s}_{m}}}\sqrt{\hat{s}_{m\vee m^*}}\,d\lambda\\&=&
H^2(\hat{s}_{m},\hat{s}_{m^*})
+2\int\pa{\sqrt{\hat{s}_{m^*}}-\sqrt{\hat{s}_{m}}}\pa{\sqrt{\hat{s}_{m\vee m^*}}
-\sqrt{\hat{s}_{m^*}}}\,d\lambda,
\end{eqnarray*}
hence, by (\ref{Eq-A1}) and Cauchy-Schwarz Inequality,
\begin{eqnarray*}
\lefteqn{H^2(\hat{s}_{m},\hat{s}_{m^*})}\hspace{7mm}\\
&\le&16[\pen(m^*)-\pen(m)]+
2\int\pa{\sqrt{\hat{s}_{m}}-\sqrt{\hat{s}_{m^*}}}\pa{\sqrt{\hat{s}_{m\vee m^*}}
-\sqrt{\hat{s}_{m^*}}}d\lambda\\&\le&16[\pen(m^*)-\pen(m)]+2H(\hat{s}_{m},\hat{s}_{m^*})
H(\hat{s}_{m\vee m^*},\hat{s}_{m^*})\\&\le&16[\pen(m^*)-\pen(m)]+\frac{1}{2}
H^2(\hat{s}_{m},\hat{s}_{m^*})+4H^2(\hat{s}_{m\vee m^*},\hat{s}_{m^*}).
\end{eqnarray*}
Therefore, for any $m\in\M$ such that $T_{m,m^*}\le0$,
\[ 
H^2(\hat{s}_{m},\hat{s}_{m^*})\le8H^2(\hat{s}_{m\vee m^*},\hat{s}_{m^*}),
+32[\pen(m^*)-\pen(m)]
\]
and, since
\begin{eqnarray*}
\lefteqn{H^2(\hat{s}_{m\vee m^*},\hat{s}_{m^*})}\hspace{7mm}\\&\leq&
4\cro{H^2(\hat{s}_{m\vee m^*},\bar{s}_{m\vee m^*}) + H^2(\bar{s}_{m\vee m^*},s)
+H^2(s,\bar{s}_{m^*})+H^2(\bar{s}_{m^*},\hat{s}_{m^*})},
\end{eqnarray*}
then
\begin{eqnarray}
(1/32)H^2(\hat{s}_{m},\hat{s}_{m^*})&\leq&H^2(\hat{s}_{m\vee m^*},\bar{s}_{m\vee
m^*}) +H^2(\hat{s}_{m^*},\bar{s}_{m^*})+\pen(m^{*})\nonumber\\&&\mbox{}
-\pen(m)+H^2(\bar{s}_{m\vee m^*},s)+H^2(s,\bar{s}_{m^*}).
\label{Eq-A0}
\end{eqnarray}
Let us set, for all $z\ge0$ and $(m,m')\in\M^2$,
\[
\Omega_z=\bigcap_{(m,m')\in\M^2}\ac{\omega\in\Omega\,\left|\,
H^2(\hat{s}_{m\vee m'},\bar{s}_{m\vee m'})\le  c|m\vee
m'|+b[\Delta_m+\Delta_{m'}+z]\right.}.
\]
It follows from (\ref{debase}) that
\begin{equation}
\P\cro{\Omega^c_z}\le
ae^{-z}\sum_{(m,m')\in\M^2}e^{-\Delta_m-\Delta_{m'}}=\Sigma^2ae^{-z}.
\label{Eq-omega}
\end{equation}
Let now $\omega$ belong to $\Omega_z$. It then follows that
\begin{equation}
H^2(\hat{s}_{m^*},\overline{s}_{m^*})\le c|m^*|+2b\Delta_{m^*}+bz\le2\pen(m^*)+bz
\label{Eq-R1}
\end{equation}
and, using Assumption~{\bf H}, that
\[ 
H^2(\hat{s}_{m\vee m^*},\bar{s}_{m\vee m^*})\le
c\delta[|m|+|m^*|]+b[\Delta_m+\Delta_{m^*}+z].
\]
Therefore we derive from (\ref{Eq-A0}), (\ref{Eq-R1}) and (\ref{Eq-pen}) that, for all $m\in\M$
such that $T_{m,m^*}\le0$,
\begin{eqnarray*}
(1/32)H^2(\hat{s}_{m},\hat{s}_{m^*})&\le&H^2(\bar{s}_{m\vee m^*},s)+H^2(s,\bar{s}_{m^*})
+(1+\delta)c|m^*|\\&&\mbox{}+3b\Delta_{m^*}+2bz+\pen(m^*)\\&\le&
H^2(\bar{s}_{m\vee m^*},s)+H^2(s,\bar{s}_{m^*})+2bz+4\pen(m^*).
\end{eqnarray*}
In order to control the bias terms $H^2(s,\bar{s}_{m'})$ of the various estimators involved in the
construction of $\tilde{s}$, we shall use Lemma~\ref{utile0} below. Since $S_{m\vee m^*}
\supset S_{m^*}$ for all $m\in\overline{\M}$, this lemma implies that
\[
H^2(\bar{s}_{m'\vee m^*},s)\le2H^2(s,S_{m'\vee m^*})\le2H^2(s,S_{m^*}),
\]
therefore
\[
H^2(\hat{s}_{m},\hat{s}_{m^*})\le128\left[H^2(s,S_{m^*})+\pen(m^*)+bz/2\right],
\]
for all $m\in\M$ such that $T_{m,m^*}\le0$ and we conclude from (\ref{Eq-A6}) and the
definition of $\D$ that, if
$\omega\in\Omega_z$,
\[
H^2(\hat{s}_{\hat m},\hat{s}_{m^*})\le\D(m^*)+\varepsilon/3\le
128\left[H^2(s,S_{m^*})+\pen(m^*)+bz/2\right]+\varepsilon/3.
\]
Since
\[
H^2(\hat{s}_{\hat m},s)\le 3\cro{H^2(\hat{s}_{\hat m},\hat{s}_{m^*})+
H^2(\hat{s}_{m^*},\overline{s}_{m^*})+H^2(\bar{s}_{m^*},s)},
\]
it follows from (\ref{Eq-R1}) and Lemma~\ref{utile0} that
\[
H^2(\hat{s}_{\hat m},s)\le3\cro{130H^2(s,S_{m^*})+130\pen(m^*)+65bz+\varepsilon/3}.
\]
Since $m^*$ is arbitrary in $\M$ we finally get 
\[
H^2(\hat{s}_{\hat m},s)\1_{\Omega_z}\leq
390\pa{\inf_{m\in\M}\cro{H^2(s,S_{m})+\pen(m)}+bz/2}+\varepsilon.
\] 
An integration with respect to $z$ taking (\ref{Eq-omega}) into account leads to
(\ref{Eq-risk1}).
\end{proof}
%
%
\begin{lem}\label{utile0} 
Within the framework of Section~\ref{A1}, for any $f\in\LL$, we have 
\[
H^2(f,\bar{f}_m)\leq 2 H^2(f,S_m)\quad\mbox{with }
\bar f_m=\sum_{I\in m\cap  \J}\left(\int_{I} f\frac{d\lambda}{\lambda(I)}\right)\1_{I}.
\]
\end{lem}
\begin{proof}
Let $\X'=\bigcup_{I\in m\cap  \J} I$. Note that $M$ is a finite measure on $\X'$ and that
for all $t\in S_m$, 
\[
H^2(f,t)=H^2(f\1_{\X'},t)+\int_{\X\setminus \X'}f d\lambda.
\]
It is therefore enough to show the result for $\X'$ in place of $\X$ and $f\1_{\X'}$ in place
of $f$ and we can restrict ourselves to the case where $M$ is a finite measure on $\X$. Let
$\sqrt{f'}$ be the $\IL_2(\X,d\lambda)$ projection of $\sqrt{f}$ on $S_m$. Since the value of
$\sqrt{f'}$ on $I$ is given by $\int_{I} \sqrt{f} d\lambda/\lambda(I)$, it suffices to prove that for each $I\in
m\cap\J$
\begin{equation}
  \label{eq:utile0}
  \int_{I} \pa{\sqrt{f}-\sqrt{\int_{I}f\frac{d\lambda}{\lambda(I)}}}^2 d\lambda\leq
2\int_{I} \pa{\sqrt{f}-\int_{I} \sqrt{f} \frac{d\lambda}{\lambda(I)}}^2d\lambda.
\end{equation}
By homogeneity, we may assume that $\lambda(I)=1$. Expanding the left-hand side
of~\eref{eq:utile0} we get
\[
\int_{I} \pa{\sqrt{f}-\sqrt{\int_{I}f d\lambda}}^2 d\lambda=
2\pa{\int_{I}fd\lambda - \int_{I}\sqrt{f}d\lambda\times \sqrt{\int_{I} f d\lambda }},
\]
which, together with the inequality $\sqrt{\int_{I} f d\lambda}\ge\int_{I} \sqrt{f} d\lambda$,
leads to the desired result.
\end{proof}
%

\section{Proofs}\label{P}

\subsection{Proof of Lemma~\ref{L-Ex4}}\label{P9}
Let $m=m_p\vee\K_j$ and $m'=m_{p'}\vee\K_{j'}$ be two elements of $\M$ and $\bar I_{p}=
\pa{[0,1)^{k}\setminus\cup_{I\in p}I}$, $\bar I_{p'}=\pa{[0,1)^{k}\setminus\cup_{I'\in p'}I'}$.
Assuming, with no loss of generality, that $j\ge j'$, we get
\[ 
m\vee m'=m_p\vee m_{p'}\vee\K_j\vee\K_{j'}=m_p\vee m_{p'}\vee\K_j
= m_1\cup m_2\cup m_3\cup m_4,
\]
with
\begin{eqnarray*}
m_1&=&\ac{K\cap I\cap I'\ne\varnothing\,|\, K\in\K_j, I\in p, I'\in p'};\\m_2&=&
\ac{K\cap I\cap\bar I_{p'}\ne\varnothing\,|\, K\in\K_j, I\in p};\\m_3&=&
\ac{K\cap\bar I_{p}\cap I'\ne\varnothing\,|\, K\in\K_j, I'\in p'};\\m_4&=&
\ac{K\cap\bar I_{p}\cap \bar I_{p'}\ne\varnothing\,|\,K\in\K_j}.
\end{eqnarray*}
Since $j<J(p)$, hence $p\subset\cup_{l>j}\K_l$, for $K\in\K_j$ and $I\in p$, $K\cap I$ is either $I$ or
$\varnothing$, so that $m=p\cup p_j$ with $p_j=\{K\cap\bar I_{p}\ne\varnothing, K\in\K_j
\}$ and
$|m|=|p|+|p_j|$. It also follows that $|m_1|\le|p|+|p'|$ and $|m_2|\le|p|$. Then,
given $K\in\K_j$ and $I'\in p'$, $K\cap I'$ is either $K$ or $I'$ or $\varnothing$ since
$K,I'\in\K$, so that $|m_3|\le|p_j|+|p'|$. Finally $|m_4|\le|p_j|$ and 
\[ 
|m\vee m'|\le2\left(|p|+|p'|+|p_j|\right)\le 2(|m|+|m'|).
\]

\subsection{Some large deviations inequalities}\label{P0}
The proofs of Theorems~\ref{densite}, \ref{poisson}, \ref{cas-discret}  and \ref{cp} require 
to check (\ref{debase}) for each specific framework. Since
\begin{equation}
H^2\left(\hat{s}_m,\bar{s}_m\right)=\sum_{I\in m\cap \J}
\pa{\sqrt{N(I)}-\sqrt{s_I}}^2\quad\mbox{for all }m\in\M,
\label{Eq-H2}
\end{equation}
this amounts to proving some deviation results for quantities of the form 
$$\sum_{I\in
m\cap\J}\pa{\sqrt{N(I)}-\sqrt{s_I}}^2-c|m|$$
 which is the purpose of this section.
Throughout it, we consider a finite set of non-negative random variables $X_I$ with $I\in
m$ and the related quantities
\begin{equation}
\chi^2(m)=\sum_{I\in m}\pa{\sqrt{X_I}-\sqrt{\E\cro{X_I}}}^2,
\label{Eq-chi2}
\end{equation}
the notation suggesting that these variables behave roughly like $\chi^2$ random variables as we
shall see. Our purpose will be to derive deviation bounds for those variables from their
expectation. Our first result is as follows:
\begin{thm}\label{chi2-I}
Let $(X_I)_{I\in m}$ be a finite set of independent non-negative random variables and
$\chi^2(m)$ be given by (\ref{Eq-chi2}). We assume that there exists $\kappa>0$ and
$\tau\geq 0$ such that 
\begin{equation}
\log\pa{\E\cro{e^{z(X_I-\E\cro{X_I})}}}\le \kappa\frac{z^2\E\cro{X_I}}{2(1-z\tau)}
\quad\mbox{for all }z\in[0, 1/\tau[,
\label{eq:lap1}
\end{equation}
and
\begin{equation}
\log\pa{\E\cro{e^{-z(X_I-\E\cro{X_I})}}}\le\kappa \frac{z^2\E\cro{X_I}}{2}
\quad\mbox{for all }z>0.
\label{eq:lap2}
\end{equation}
Let 
\[
K=\max\ac{\sqrt{2}\ ;\ \frac{\sqrt{2}}{2} +\sqrt{\pa{\frac{\tau}{\kappa}-\frac{1}{2}}_+}}.
\]
Then for all $x>0$,
\begin{equation}\label{cle}
\P\cro{\chi^2(m)\ge\E\cro{\chi^2(m)}+K^2\kappa\pa{2\sqrt{2|m|x}+ x}}\le e^{-x},
\end{equation}
and
\begin{equation}\label{cle2}
\P\cro{\chi^2(m)\le\E\cro{\chi^2(m)}-2K^2\kappa\sqrt{2|m|x}}\leq e^{-x}.
\end{equation}
\end{thm}

\begin{proof}
Let us first introduce the following large deviation result, the proof of which follows the
lines of the proof of Lemma~8 of Birg\'e and Massart (1998).
%
\begin{lem}\label{lemLD}
Let $Y_1,\ldots,Y_n$ be $n$ independent, centered random variables. If
\[
\log\left(\E\cro{e^{z Y_i}}\right)\le\kappa\frac{z^2\theta_i}{2(1-z \tau)}
\quad\mbox{for all }z\in[0, 1/\tau[\quad\mbox{and}\quad1\le i\le n,
\]
then
\[
\P\cro{\sum_{i=1}^nY_i\ge\left(2\kappa
    x\sum_{i=1}^n\theta_i\right)^{1/2} + \tau x}
\leq e^{-x}\quad\mbox{for all }x>0.
\]
If, for $1\le i\le n$ and all $z>0$, $\log\left(\E\cro{e^{-z Y_i}}\right)\le\kappa
z^2\theta_i/2$, then
\[
\P\cro{\sum_{i=1}^n Y_i\le -\left(2\kappa x\sum_{i=1}^n\theta_i\right)^{1/2}}\leq 
e^{-x}\quad\mbox{for all }x>0.
\]
\end{lem}
It follows from~\eref{eq:lap1}, \eref{eq:lap2} and Lemma~\ref{lemLD} with $n=1$,
$Y_1=X_I-\E\cro{X_I}$ and $\theta_1= \E\cro{X_I}$ that, for all $x>0$ and $I\in m$,
\[
\P\cro{X_I\geq \E\cro{X_I} +\sqrt{2\kappa \E\cro{X_I} x} + \tau x}\leq e^{-x}
\]
and
\[
\P\cro{X_I\leq \E\cro{X_I} -\sqrt{2\kappa \E\cro{X_I} x}}\leq e^{-x}.
\]
Setting $u=\E\cro{X_I}/(\kappa x)$, we deduce that, with probability not smaller than
$1-2e^{-x}$,
\begin{eqnarray*}
\lefteqn{\ab{\sqrt{X_I}-\sqrt{\E\cro{X_I}}}}\hspace{15mm}\\&\leq&
\max\left\{\sqrt{\E\cro{X_I}}-
\sqrt{\pa{\E\cro{X_I} -\sqrt{2\kappa\E\cro{X_I}x}}_{+}};\right.\\&&\mbox{}
\hspace{10mm}\left.\sqrt{\E\cro{X_I}+\sqrt{2\kappa\E\cro{X_I} x} +\tau x}-
\sqrt{\E\cro{X_I}}\right\}\\&=&\sqrt{\kappa x}\,\max\ac{\sqrt{u}-
\sqrt{\pa{u-\sqrt{2u}}_{+}};\,\sqrt{u+\sqrt{2u}+(\tau/\kappa)}-\sqrt{u}}\\&\le&
\sqrt{\kappa x}\,\sup_{z>0}\max\ac{\sqrt{z}-\sqrt{\pa{z -\sqrt{2z}}_{+}};\ 
\sqrt{z+\sqrt{2z}+(\tau/\kappa)}-\sqrt{z}}.
\end{eqnarray*}
On the one hand, note that $z\to \sqrt{z}-\sqrt{\pa{z -\sqrt{2z}}_{+}}$ admits a
maximum equal to $\sqrt{2}$ for $z=2$. On the other hand, using the
inequality $\sqrt{a+b}\leq \sqrt{a}+\sqrt{b}$ which holds for all
positive numbers $a,b$, we obtain for all $z>0$, 
\begin{eqnarray*}
 \sqrt{z+\sqrt{2z}+(\tau/\kappa)}-\sqrt{z}&\leq&
 \sqrt{\pa{\sqrt{z}+\frac{\sqrt{2}}{2}}^2+\pa{\frac{\tau}{\kappa}-\frac{1}{2}}_+}-
\sqrt{z}\\ &\leq& \frac{\sqrt{2}}{2} + \sqrt{\pa{\frac{\tau}{\kappa}-\frac{1}{2}}_+}
\end{eqnarray*}
and therefore $\ab{\sqrt{X_I}-\sqrt{\E\cro{X_I}}}\leq K\sqrt{\kappa x}$ with
probability not smaller than $1-2e^{-x}$, or equivalently
\begin{equation}\label{interm1}
\P\cro{U_I\geq K^2x}\leq 2e^{-x}\quad\mbox{for all }x>0\quad\mbox{with}\quad
U_I=\kappa^{-1}\pa{\sqrt{X_I}-\sqrt{\E\cro{X_I}}}^2.
\end{equation}
Since $\chi^2(m)=\kappa\sum_{I\in m}U_I$ and the random variables $U_I,I\in m$ are
independent,  (\ref{cle}) will derive from Lemma~\ref{lemLD} if we show, setting
$E_I=\E\cro{U_I}$, that 
\begin{equation}\label{laplace1}
\log\pa{\E\cro{e^{z(U_I-E_I)}}}\leq \frac{4K^4 z^2}{2(1-K^2z)}
\quad\mbox{for all }z\in]0,1/K^2[.
\end{equation}
Similarly, (\ref{cle2}) will follow from
\begin{equation}\label{laplace2}
\log\pa{\E\cro{e^{-z(U_I-E_I)}}}\leq \frac{4K^4 z^2}{2}\quad\mbox{for all }z>0.
\end{equation}

To prove~\eref{laplace1}, we shall use the following lemma about the centered moments
of  positive random variables.
\begin{lem}
Let $Z$ be a non-negative random variable. For any positive even integer $k$,
$$\E\cro{\pa{Z-\E\cro{Z}}^k}\leq \E\cro{Z^k}-(\E\cro{Z})^k\leq \E\cro{Z^k}.$$
\end{lem}
Note that the inequality $\E\cro{\pa{Z-\E\cro{Z}}^k}\leq\E\cro{Z^k}$ also holds 
true for odd integers $k$ since $\E\cro{Z}\geq 0$ and the map $z\mapsto z^k$ is
then increasing.
\begin{proof}
Since the result is trivial for $k=2$, we may assume that $k\ge4$ and, using homogeneity,
that $\E\cro{Z}=1$. Consider the function $z\mapsto  Q(z)= z^k-(z-1)^k-k(z-1)$ on
$[0,+\infty[$. Its second derivative is negative for $z<1/2$ and positive for $z>1/2$, from
which we easily derive that $Q$ has a minimum for $z=1$. This shows that $Q(z)\ge 1$ for
all $z\ge 0$ and consequently, 
\[
\E\cro{Z^k}-\E\cro{\pa{Z-1}^k}=\E\cro{Q\pa{Z}}\geq Q(1)=1
\]
which leads to the result.
\end{proof}
The random variable $U_I$ is positive and by~\eref{interm1} satisfies
$\P\cro{U_I\ge t}\leq 2e^{-t/K^2}$. Consequently, we deduce from 
the previous lemma (with $Z=U_I$) that for all integers $k$ (odd or even)
\begin{equation}
\E\cro{\pa{U_I-E_I}^k}\le\E\cro{U_I^k}=
\int_{0}^{+\infty}kt^{k-1}\P\cro{U_I\ge t}dt\le2(k!)K^{2k}.
\label{Eq-A5}
\end{equation}
Hence, for all $z\in]0,1/K^2[$\,,
\[
\log\left(\E\cro{e^{z(U_I-E_I)}}\right)\le
\log\pa{1+0+2\sum_{k\ge2}z^kK^{2k}}\le 2\sum_{k\ge2}z^kK^{2k}=
\frac{4K^4 z^2}{2(1-K^2z)}.
\]
To prove ~\eref{laplace2}, note that, for all $z,u>0$, $e^{-z u}\leq 1-z u + z^2u^2/2$.
Therefore, by (\ref{Eq-A5}),
\[
\log\left(\E\cro{e^{-z(U_I-E_I)}}\right)=\log\left(\E\left[e^{-zU_I}\right]\right)
+zE_I\le\frac{z^2}{2}\E\cro{U_I^2}\leq \frac{4K^4z^2}{2},
\]
which completes the proof of Theorem~\ref{chi2-I}.
\end{proof}

A second pair of deviation inequalities for variables of the form $\chi^2(m)$ is as follows.
\begin{thm}\label{chi2-II}
Let $m$ be a finite index set and $\bm{X}_j=\pa{X_{I,j}}_{I\in m}$, $1\le j\le p$ be i.i.d.\ random
vectors with values in $\R_+^{|m|}$. Assume that there exist positive numbers $A$ and
$\kappa$ such that
\begin{equation}\label{eq:var-ct0}
\sum_{I\in m} X_{I,1}\leq A\;\; a.s.\qquad\mbox{and}\qquad
{\rm Var}\pa{X_{I,1}}\le\kappa\E\cro{X_{I,1}}\quad\mbox{for all }I\in m.
\end{equation}
If $X_I=\sum_{j=1}^pX_{I,j}$ for all $I\in m$ and $\chi^2(m)$ is given by (\ref{Eq-chi2}), then
\begin{equation}\label{eq:chi2}
\P\cro{\chi^2(m)\geq 8\kappa |m|+ 202Ax}\leq e^{-x}
\quad\mbox{for all }x>0.
\end{equation}
\end{thm}
\begin{proof}
Since $X_{I,1}=0$ a.s.\ if $\E[X_{I,1}]=0$, we may remove all indexes $I$ such that $\E[X_{I,1}]=0$
in the sum and therefore assume that $\E\cro{X_I}=p\E\cro{X_{I,1}}>0$ for all $I\in m$. We can
then write, for all $z>0$, 
\begin{eqnarray*}
\lefteqn{\P\pa{\sqrt{\X^2(m)}\geq z}}\\&=& \P\pa{\sum_{I\in m}
\frac{\pa{\sqrt{X_I}-\sqrt{\E\cro{X_I}}}}{\sqrt{\X^2(m)}}
\pa{\sqrt{X_I}-\sqrt{\E\cro{X_I}}}\ge z,\sqrt{\X^2(m)}\geq z}\\&=& 
\P\pa{\sum_{I\in m}\frac{\pa{\sqrt{X_I}-\sqrt{\E\cro{X_I}}}}{\sqrt{\X^2(m)}}
\frac{X_I-\E\cro{X_I}}{\sqrt{X_I}+\sqrt{\E\cro{X_I}}}\ge z,\sqrt{\X^2(m)}\geq z}\\
&=& \P\pa{\sum_{j=1}^p \!\!\cro{\sum_{I\in m}\frac{
\pa{\sqrt{X_I}-\sqrt{\E\cro{X_I}}}\pa{X_{I,j}-\E\cro{X_{I,j}}}}
{\sqrt{\X^2(m)}\pa{\sqrt{X_I}+\sqrt{\E\cro{X_I}}}}}
\ge z,\sqrt{\X^2(m)}\ge z}\\&=& \P\pa{\sum_{I\in m}\!\!
\cro{\sum_{j=1}^p \frac{\pa{\sqrt{X_I}-\sqrt{\E\cro{X_I}}}\sqrt{\E\cro{X_I}}}
{\sqrt{\X^2(m)}\pa{\sqrt{X_I}+\sqrt{\E\cro{X_I}}}}
\frac{X_{I,j}-\E\cro{X_{I,j}}}{\sqrt{\E\cro{X_I}}}}\!\!\ge z,\sqrt{\X^2(m)}\geq z}\\
&=& \P\pa{\sum_{j=1}^p \sum_{I\in m}
t_{I}\frac{X_{I,j}-\E\cro{X_{I,j}}}{\sqrt{\E\cro{X_I}}}\ge z, \sqrt{\X^2(m)}\geq z},
\end{eqnarray*}
where 
\[
t_{I}=\frac{\pa{\sqrt{X_I}-\sqrt{\E\cro{X_I}}}\sqrt{\E\cro{X_I}}}
{\sqrt{\X^2(m)}\pa{\sqrt{X_I}+\sqrt{\E\cro{X_I}}}}\qquad\mbox{for all }I\in m .
\]
Note that $\sum_{I\in m} t_I^2\leq 1$ since $\sqrt{\E\cro{X_I}}/(\sqrt{X_I}+
\sqrt{\E\cro{X_I}})\leq 1$ and that $|t_I|\leq  z^{-1}\sqrt{\E\cro{X_I}}$ on the set
$\sqrt{\X^2(m)}\geq z$, from which we deduce that
\begin{equation}
\P\pa{\sqrt{\X^2(m)}\geq z}\leq \P\pa{\sup_{\bm{t}\in \T} \sum_{j=1}^p \sum_{I\in m}
t_{I}\frac{X_{I,j}-\E\cro{X_{I,j}}}{\sqrt{\E\cro{X_I}}}\ge z},
\label{Eq-A7}
\end{equation}
where $\T$ denotes the set of vectors $\bm{t}=\pa{t_I}_{I\in m}\in\R^{|m|}$
satisfying 
\begin{equation}
|t_I|\leq \frac{\sqrt{\E\cro{X_I}}}{z}\quad\mbox{for all }I\in m
\qquad\mbox{and}\qquad\sum_{I\in m} t_I^2\leq 1.
\label{Eq-A8}
\end{equation}
In order to bound the right-hand side of (\ref{Eq-A7}), we shall use the following result
from Massart (2000, Theorem~2.4).
\begin{thm}\label{massart}
Let $\bm{\xi}_1,\ldots,\bm{\xi}_p$ be independent random variables with values in some
measurable space $\H$ and $\F$ be some countable family of real valued measurable
functions on $\H$ such that $\|f\|_{\infty}\le b<+\infty$ for all $f\in\F$.  If
\[
Z=\sup_{f\in\F}\ab{\sum_{j=1}^p
f(\bm{\xi}_j)-\E\cro{f(\bm{\xi}_j)}}\qquad\mbox{and}\qquad
\sigma^2=\sup_{f\in\F}\cro{\sum_{j=1}^p {\rm Var}\pa{f(\bm{\xi}_j)}},
\]
then for every positive numbers $\eps ,x$
\[
\P\cro{Z\geq (1+\eps) \E\cro{Z}+ 2\sigma\sqrt{2x} +
\left(2.5+32\varepsilon^{-1}\right)b x}\leq e^{-x}.
\]
\end{thm}
We want to apply this result to the vectors $\bm{\xi}_j\in\R^{|m|}$ with coordinates
$\xi_{I,j}=(X_{I,j}-\E\cro{X_{I,j}})/\sqrt{\E\cro{X_I}}$ for $I\in m$. Under our
assumptions, these random vectors are independent and satisfy 
\[
\sum_{I\in m}\sqrt{\E\cro{X_I}}|\xi_{I,j}|\leq \sum_{I\in m} \pa{X_{I,j}
+\E\cro{X_{I,j}}}\leq 2A.
\] 
Consequently, the random vectors $\bm{\xi}_j$ take their values in the subset $\H$ of 
$\R^{|m|}$ given by 
\[
\H=\ac{\bm{u}=(u_I,\ I\in m)\;\left|\;\sum_{I\in m}\sqrt{\E\cro{X_I}}|u_I|\le2A\right.}.
\]
For $\bm{u}\in\H$ and $\bm{t}\in\T$, we set $f_{\bm{t}}(\bm{u})=\sum_{I\in m}t_Iu_I$ and 
$\F=\ac{f_{\bm{t}}, \bm{t}\in\T'}$ where $\T'$ denotes a countable and dense subset of $\T$.
With no loss of generality we can assume that $\T'$ is symmetric around 0 (if $\bm{t}\in\T'$ then
$-\bm{t}\in\T'$) which implies that the absolute values can be removed in the definition of $Z$.
Since, for all $\bm{t}\in\T$ and $1\le j\le p$, $f_{\bm{t}}(\bm{\xi}_j)$ is centered, we can finally write
\[
Z=\sup_{\bm{t}\in \T} \sum_{j=1}^p \sum_{I\in m}t_I\frac{X_{I,j}-\E\cro{X_{I,j}}}
{\sqrt{\E\cro{X_I}}}=\sup_{\bm{t}\in \T}\sum_{I\in m}t_I
\pa{\sum_{j=1}^p\frac{X_{I,j}-\E\cro{X_{I,j}}}{\sqrt{\E\cro{X_I}}}}.
\]
Using Cauchy-Schwarz Inequality and (\ref{Eq-A8}), we then derive that
\[
\E^2\cro{Z}\le\E\cro{Z^2}\le\sum_{I\in m}\E\cro{\pa{\sum_{j=1}^p
  \frac{X_{I,j}-\E\cro{X_{I,j}}}{\sqrt{\E\cro{X_I}}}}^2}
=\sum_{I\in m} \sum_{j=1}^p \frac{{\rm Var}(X_{I,j})}{\E\cro{X_I}}.
\]
Since ${\rm Var}(X_{I,j})\le\kappa \E\cro{X_{I,j}}$ and $\sum_{j=1}^p\E\cro{X_{I,j}}
=\E\cro{X_I}$, we conclude that $\E\cro{Z}\leq \sqrt{\kappa|m|}$.
To bound $\|f_{\bm{t}}\|_\infty$, we use (\ref{Eq-A8}) which implies that, for all $\bm{u}\in
\H$ and $\bm{t}\in\T$,
\[
|f_{\bm{t}}(\bm{u})|=\left|\sum_{I\in m}t_{I}u_{I}\right|\leq \sum_{I\in m}
|t_{I}||u_{I}|\leq \sum_{I\in m}\frac{\sqrt{\E\cro{X_I}}|u_I|}{z}\leq \frac{2A}{z}.
\]
Finally,  it follows from the equidistribution of the $\bm{X}_j$, Cauchy-Schwarz Inequality, 
(\ref{eq:var-ct0}) and (\ref{Eq-A8}) that, for all $\bm{t}\in \T$, 
\begin{eqnarray*}
\sum_{j=1}^p {\rm Var}\pa{f_{\bm{t}}(\bm{\xi}_j)}&=&p{\rm Var}\pa{f_{\bm{t}}(\bm{\xi}_1)}
\;\;=\;\;p\E\cro{\pa{\sum_{I\in m}t_I\frac{X_{I,1}-\E\cro{X_{I,1}}}
{\sqrt{p\E\cro{X_{I,1}}}}}^2}\\&\le&2 \E\cro{\pa{\sum_{I\in m}t_{I}\frac{X_{I,1}}
{\sqrt{\E[X_{I,1}]}}}^2}+ 2\pa{\sum_{I\in m}t_I\sqrt{\E[X_{I,1}]}}^2\\&\le&
2\E\left[\left(\sum_{I\in m}X_{I,1}\right)\left(\sum_{I\in m}t_I^2\frac{X_{I,1}}
{\E[X_{I,1}]}\right)\right]+2\sum_{I\in m}t_{I}^2\sum_{I\in m}\E[X_{I,1}]\\&\le&2A
\left(\E\left[\sum_{I\in m}t_I^2\frac{X_{I,1}}{\E[X_{I,1}]}\right]+\sum_{I\in m}t_I^2\right)
\;\;\le\;\;4A.
\end{eqnarray*}
In view of all these bounds, we may apply Theorem~\ref{massart} with $\sigma^2=4A$,
$b=2A/z$ and $\eps=1$ and obtain that $\P\cro{\sqrt{\chi^2(m)}\geq z}\leq e^{-x}$ as soon 
as $z\geq 2\sqrt{\kappa|m|}+4\sqrt{2Ax}+69Ax/z$.  Solving this quadratic inequation and
using  $(a+b)^2\le2\left(a^2+b^2\right)$, we can check that this inequality holds if
$z^2\ge8\kappa|m|+202Ax$, hence the result.
\end{proof}
\subsection{Density estimation}\label{P1}
%
\subsubsection{Proof of Theorem~\ref{densite}\label{P1a}}
For two given classes $m,m'\in\M$, we apply Theorem~\ref{chi2-II} with  $m''=m\vee
m'$ in place of $m$, $p=n$ and $X_{I,j}=\1_{Y_j\in I}$ for all $I\in m''$ and $j=1,\ldots,n$.
Then $X_I=nN(I)$ and (\ref{eq:var-ct0}) is satisfied with  $A=\kappa=1$ since $X_{I,1}$  is
a Bernoulli random variable and we derive from (\ref{Eq-H2}) that, for all $x>0$, with
probability not smaller than $1-e^{-x}$,
\[
H^2(\hat s_{m''},\overline{s}_{m''})=\sum_{I\in m''}\pa{\sqrt{N(I)}-\sqrt{\E\cro{N(I)}}}^2=
\frac{\chi^2(m'')}{n}\le \frac{8|m''|+202x}{n}.
\]
Therefore \eref{debase} holds with $c=8/n$, $a=1$ and $b=202/n$. We then conclude from
Theorem~\ref{main} and the fact that $H^2(t,u)$ is always bounded by 2.
%
\subsubsection{Proof of Proposition~\ref{P-ex3}\label{P1b}}
By assumption, $\sqrt{s}$ has a variation bounded by $R$ and we may apply to it
Corollary~1 of Barron, Birg\'e and Massart (1999) with $\alpha=1$, $D=2^j$ with $j\ge2$
and $N=2^{3j}$. It follows that one can find $m\in\M_{3j,D}$ such that
$H^2(s,S_m)\le(64/3)(R/D)^2$. Since $\pen(m)\le CjDn^{-1}$ for $m\in\M_{3j,D}$, we derive
from Theorem~\ref{densite} that
\[
\E_s\left[H^2(\tilde{s},s)\right]\le C'\inf_{j\ge2}\left\{R^22^{-2j}+j2^jn^{-1}\right\}.
\]
Then (\ref{Ri2}) follows if we define $j\ge2$ by
\[
4^{-j+1}\le\left[nR^2/\log\left(1+nR^2\right)\right]^{-2/3}<4^{-j+2},
\]
which is always possible since $nR^2>0$, and distinguish between the cases $j=2$
(which corresponds to $nR^2\le26.519$) and $j>2$.

When $\sqrt{s}$ is continuous with modulus $w$, there exists an element $t\in S_{m_j}$
such that $\|\sqrt{s}-\sqrt{t}\|_\infty\le w(2^{-j})$, hence  $H(s,S_{m_j})\le w(2^{-j})$. Since
$x_w>0$, we can choose $j$ such that $2^{-j}<x_w\le2^{-j+1}$. Recalling that $\pen(m_j)\le
C2^j/n$, we deduce from Theorem~\ref{densite} that
\[
\E_s\left[H^2(\tilde{s},s)\right]\le C'\left[w^2(2^{-j})+n^{-1}2^j\right]\le C'
\left[w^2(x_w)+2(nx_w)^{-1}\right]\le 3C'(nx_w)^{-1},
\]
which proves (\ref{Ri1}). If $\sqrt{s}$ belongs to $\H_\alpha^R$ with $R\ge n^{-1/2}$,
then $x_w=(nR)^{-2/(2\alpha+1)}$ and the risk bound follows.

If $s$ belongs to $\S^3(D,R)$, we can write  $s=\sum_{k=1}^D s_k\1_{[x_{k-1},x_k)}$ with
$0=x_0<x_1<\ldots<x_D=1$ and $\sup_{1\le k\le D}s_k\le R$. Fix $l$ such that
$2^l\ge nR>2^{l-1}$. Then $2^l\ge2D$ and for $0\le k\le D$, set $x'_k=
\sup\{x\in\J_l\,|\,x\le x_k\}$ and $t=\sum_{k=1}^D s_k\1_{[x'_{k-1},x'_k)}$ so that $t\in
S_m$ with $m\in\M_{l,D'}$ with $D'\le D$ since some intervals $[x'_{k-1},x'_k)$ may be empty.
Then
\[
H^2(s,t)\le R\sum_{k=1}^{D-1}(x_k-x'_k)<RD2^{-l}.
\]
Recalling from (\ref{Eq-Dm4}) that $\pen(m)\le Cn^{-1}[D(l\log 2+2-\log D)+2\log l]$ for
$m\in\M_{l,D}$, we conclude from Theorem~\ref{densite}, (\ref{Eq-Dm4}) and our
choice of $l$ that
\begin{eqnarray*}
\lefteqn{\E_s\left[H^2(\tilde{s},s)\right]}\hspace{15mm}\\
&\le&C'\left[RD2^{-l}+[D(l\log 2+2-\log D)+2\log l]n^{-1}
\right]\\&\le&C'(D/n)\left[3+\log2+\log\left(2^{l-1}/D\right)+2D^{-1}\log l\right]
\\&\le&C'(D/n)\left[3+\log2+\log\left(nR/D\right)+2(D\log2)^{-1}\log\log(2nR)\right]
\end{eqnarray*}
and (\ref{Ri3}) follows since $nR\ge 2D$.

\subsection{Random vectors}\label{P3}
%
\subsubsection{Proof of Theorem~\ref{cas-discret}\label{P3a}}
For two given elements $m,m'\in\M$, we apply Theorem~\ref{chi2-I} with
$m''= m\vee m'$ in place of $m$ and $X_I=N(I)$. We derive from the independence of the
$N_i$ that~\eref{eq:lap1} and~\eref{eq:lap2} hold. Therefore, for all $x>0$, with
probability not smaller than $1-e^{-x}$,
\begin{eqnarray*}
H^2(\hat s_{m''},\overline{s}_{m''})&=&\sum_{I\in m''}
\pa{\sqrt{N(I)}-\sqrt{\E\cro{N(I)}}}^2\\ &\leq& \E\cro{\sum_{I\in m''}
\pa{\sqrt{N(I)}-\sqrt{\E\cro{N(I)}}}^2} + K^2\kappa\pa{2\sqrt{2|m''|x}+x}.
\end{eqnarray*}
If follows from~\eref{eq:lap1} that ${\rm Var}(N(I))\leq \kappa\E\cro{N(I)}$ (expand
both side of ~\eref{eq:lap1} in a vicinity of 0) and therefore
\begin{eqnarray*}
\E\cro{\sum_{I\in m''}\pa{\sqrt{N(I)}-\sqrt{\E\cro{N(I)}}}^2}&=&\sum_{I\in m''}
\E\cro{ \pa{\sqrt{N(I)}-\sqrt{\E\cro{N(I)}}}^2}\\ &\leq& \sum_{I\in m''}
\E\cro{\frac{\pa{N(I)-\E\cro{N(I)}}^2}{\E\cro{N(I)}}}\;\;\le\;\;\kappa |m''|.
\end{eqnarray*}
Using the inequality $2\sqrt{2|m''|x}\leq |m''|+2x$ we conclude that, with probability
not smaller than $1-e^{-x}$,
\begin{equation}
H^2(\hat s_{m''},\overline{s}_{m''})\le\pa{1+K^2}\kappa|m''|+3K^2\kappa x.
\label{Eq-H9}
\end{equation}
We derive that~\eref{debase} is fulfilled with $c=\pa{1+K^2}\kappa$, $b=3K^2\kappa$,
$a=1$ and Theorem~\ref{cas-discret} follows from Theorem~\ref{main}.

%
\subsubsection{Proof of Proposition~\ref{P-mon6}\label{P3b}}
Let us first note that, if $|m|=n$, then $H^2(s,S_m)=0$, hence by (\ref{eq-78}),
$\E\left[H^2(\tilde{s},s)\right]\le C(\kappa,K)n$ which proves the bound when
$R>n/\sqrt{3}$. For the other cases, we deduce from Lemma~\ref{monotone} below that, for
any  $D\in\X$, one can find some $m\in\M$ such that $|m|\le D$ and $H^2(s,S_m)\le
n(R/D)^2$. Setting $D=1$, we get the result for the case $R^2<n^{-1}\log n$. Finally,
when $n^{-1}\log n\le R^2\le n^2/3$ we fix $D=\inf\left\{j\in\Bbb{N}\,|\,j^3\ge
nR^2/\log(n/R)\right\}$. Since the function $R\mapsto R^2/\log(n/R)$ is increasing for
$R<n/\sqrt{3}$, $1\le D\le n$ and the corresponding risk bound follows.

\begin{lem}\label{monotone}
Let $f$ be a nondecreasing function from $\X=\{1,\ldots,n\}$ to $\R$ such that 
$\sqrt{f(n)}-\sqrt{f(1)}=R$. For $D\in\X$, one can find a partition $(I_1,\ldots,I_K)$ of $\X$ into
$K\le D$ intervals and a function $g$ from $\X$ to $\R$ of the form
$g=\sum_{k=1}^K\beta_k\1_{I_k}$  such that
\[ 
\sum_{i=1}^n\left(\sqrt{f(i)}-\sqrt{g(i)}\right)^2\le nR^2 D^{-2}.
\] 
\end{lem}
%
\noindent{\em Proof:}
Let us set $j_0=1$ and define iteratively for $k\ge1$, using the convention $\inf\varnothing =n$,
\begin{equation}
j_{k}=\inf\ac{j\in\ac{j_{k-1}+1,\ldots,n}\,\left|\,\sqrt{f(j)}-\sqrt{f(j_{k-1})}>R/D\right.}.
\label{Eq-Jk}
\end{equation}
Let $K=\inf\ac{k\geq 1,\ j_{k}=n}$, $I_K=\ac{j_{K-1},\ldots,n}$ and  for
$k=1,\ldots,K-1$ (if $K\geq 2$), $I_{k}=\ac{j_{k-1},\ldots,j_{k}-1}$. This defines a partition
of $\X$ with $K$ elements and it follows from (\ref{Eq-Jk}) that
\[
R=\sqrt{f(n)}-\sqrt{f(1)}\ge\sum_{k=1}^{K-1}\sqrt{f(j_k)}-\sqrt{f(j_{k-1})}> (K-1)R/D,
\] 
hence $K-1<D$ and $K\le D$. Let us now set $\beta_k=f(j_{k-1})$ for $1\le k\le K$. Since
$\sqrt{f(j_k-1)}-\sqrt{f(j_{k-1})}\le R/D$ we get for all $i\in I_k$, $0\le\sqrt{f(i)}-\sqrt{g(i)}\le
R/D$. Hence,
\[
\sum_{i=1}^n\left(\sqrt{f(i)}-\sqrt{g(i)}\right)^2=\sum_{k=1}^{K}\sum_{i\in I_k}
\pa{\sqrt{f(i)}-\sqrt{g(i)}}^2\le nR^2D^{-2}.\cqfd
\]
%

\subsection{Poisson and other counting processes}\label{P4}
%
\subsubsection{Poisson processes\label{P4t}}
The proof of Theorem~\ref{poisson} follows the same lines as the proof of
Theorem~\ref{cas-discret}. We apply Theorem~\ref{chi2-I} with $m''= m\vee m'$ in place of $m$
and $X_I=N(I)$. Since $\ac{N(I),\ I\in m''}$ are independent Poisson random variables, the
assumptions of the theorem are fulfilled with $\kappa=\tau=1$. We then proceed as for
Theorem~\ref{cas-discret} to get (\ref{Eq-H9}) with $K^2=2$ which provides the relevant values
of $c$ and $b$. 

\subsubsection{Proof of Theorem~\ref{cp}\label{P4a}}
Let us fix two classes $m,m'\in\M$. We first apply Theorem~\ref{chi2-II} with $m''=m\vee m'$ in
place of $m$, $p=n$ and $X_{I,j}=N^j(I)$ for all $I\in m''$ and $j=1,\ldots,n$. Then for all $I\in
m''$, $N(I)=X_{I}$. Since $X_{I,j}$ is bounded by $k$, $\E\cro{X_{I,j}^2}\le k\E\cro{X_{I,j}}$ and 
(\ref{eq:var-ct0}) holds with $A=\kappa=k$. This implies that, for all $x>0$, with probability not
smaller than $1-e^{-x}$, 
\begin{equation}\label{eq:etape1}
\sum_{I\in m''}\pa{\sqrt{N(I)}-\sqrt{\E\cro{N(I)}}}^2 \le k\pa{8|m''| + 202x}.
\end{equation}
Then we apply once again Theorem~\ref{chi2-II} with $m''=m\vee m'$ in place of $m$, $p=n$ and
$X_{I,j}=\int_I sY^jd\lambda$ for all $I\in m''$ and $j=1,\ldots,n$. Since $Y^j$ is bounded
by 1, the assumptions of
Theorem~\ref{chi2-II} are fulfilled with $A=\int_{\X} s\,d\lambda$ and $\kappa=\kappa'$. 
Consequently, with probability not smaller than $1-e^{-x}$,
\begin{equation}
\label{eq:etape2}
\sum_{I\in m''}\pa{\sqrt{\int_{I}sYd\lambda}-\sqrt{\E\cro{\int_{I}sd\lambda}}}^2\leq
8\kappa'|m''| + 202Ax. 
\end{equation}
Since $\E\cro{\int_{I}sd\lambda}=\E\cro{N(I)}$, we derive from \eref{eq:etape1} and \eref{eq:etape2}
that, with probability not smaller than $1-2e^{-x}$,
\begin{eqnarray*}
\lefteqn{H^2(\hat s_{m''},\overline{s}_{m''})}\hspace{5mm}\\&\leq& \sum_{I\in m''}
\pa{\sqrt{N(I)}-\sqrt{\int_{I}sYd\lambda}}^2\\
&\leq& 2\sum_{I\in m''}\pa{\sqrt{N(I)}-\sqrt{\E\cro{N(I)}}}^2 + 2\sum_{I\in m''}
\pa{\sqrt{\int_{I}sYd\lambda}-\sqrt{\E\cro{\int_{I}sd\lambda}}}^2\\
&\le& 16|m''|(k+\kappa')+404x(k+A).
\end{eqnarray*}
This means that \eref{debase} holds  with $c=16(k+\kappa')$, $a=2$ and $b=404(k+A)$.
Therefore, if we set $\Delta_m=k(k+A)^{-1}\Delta'_m$ for all $m\in\M$, (\ref{som})
holds with $\Sigma=\Sigma'(k/(k+A))$ and
$\pen(m)=16\delta|m|(k+\kappa')+404k\Delta'_m$.  An application of
Theorem~\ref{main} leads to the result.

\subsubsection{Proof of Proposition~\ref{P-surv}\label{P4b}}
The following argument shows that (\ref{Eq-fond1}) is satisfied: let $A$ be some measurable
subset of $\X$ and $B$ be the subset of $A$ given by $B=\ac{t\in A\  |\
\lambda\pa{[0,t]\cap A}=0 }$. Since, by definition, the sets $[0,t]\cap B$ with $t\in B$ are
negligible,  $\lambda(B)=0$ (write $B$ as an at most countable union of those sets).
Consequently, 
\begin{eqnarray*}
\P\pa{N(A)>0,\ M(A)=0}
&\leq&\sum_{j=1}^n \P\pa{N^j(A)=1,\ \int_{A}\1_{\widetilde{T}_j\geq t} dt=0}\\&\leq&
\sum_{j=1}^n \P\pa{\widetilde{T}_j=T_j,\ T_j\in A,\ \lambda\pa{A\cap[0,\widetilde{T}_j]}=0}\\&\le&
\sum_{j=1}^n \P\pa{T_j\in B}\;\;=\;\;0
\end{eqnarray*}
since the common distribution of the $T_j$ is continuous. Moreover
\[
\int_{\X} sd\lambda=\int_0^1\frac{p(t)}{\Bbb{P}[T_1\ge t]}dt=-\log([\P(T_1\geq 1)])
\]
since $-p(t)$ is the derivative of $\Bbb{P}[T_1\ge t]$. Finally we can take $\kappa'=2$ since,
whatever $I\subset\X$,
\begin{eqnarray*}
{\rm Var}\cro{\int_I s(t)Y_t^*dt}&\le& \E\cro{\pa{\int_Is(t)Y_t^*dt}^2}\;\;=\;\;
\E\cro{\int_{I\times I}s(t)s(t')Y_t^*Y_{t'}^*dt\,dt'}\\
&=& \int_{I\times I}s(t)s(t')\E\cro{Y_t^*Y_{t'}^*} dt\,dt'\\
&=& \int_{I\times I}s(t)s(t')\,\P\left[\widetilde{T}_1\geq
\max\ac{t,t'}\right]\,dt\,dt'\\&=&  2\int_I s(t)\pa{\int_I\1_{\{t'\ge t\}}\,s(t')\,
\P\left[\widetilde{T}_1\geq t'\right]dt'}dt\\&\le& 2 \int_Is(t)\,\E\cro{\int_t^1
s(t')Y_{t'}^*\,dt'}dt\\ &=& 2 \int_Is(t) \E\cro{N^1([t,1])}dt\\&\leq& 2 
\int_I s(t) \P\left[\widetilde{T}_1\geq t\right]\,dt\;\;=\;\; 2\E\cro{\int_I s(t)Y_t^* dt}.
\end{eqnarray*}
%
\subsubsection{Proof of Proposition~\ref{P-mark}\label{P4c}}
Clearly (\ref{aalen}) holds true. We now prove that  Condition \eref{Eq-fond1} is also
fulfilled. Let $A$ be some measurable subset of $\R_+$ and for $l\geq 1$ let $B_l$
be the subset of $A$ defined by 
\[
B_l=\ac{t\in A\ |\ \lambda\pa{]t-l^{-1},  t]\cap A}=0}.
\]
For each $l\geq 1$, note that the sets $[t-l^{-1},t]\cap B_{l}\subset [t-l^{-1},t]\cap A$ are
negligible for $t\in B_l$ and hence so is $B_l$ (write $B_l$ as an at most countable union of
those). Denoting, for $j=1,\ldots,n$, the time of the jump of $X^j$ from state 1 to 0  by
$T^j_{1,0}$, we have
\begin{eqnarray*}
\lefteqn{\P\pa{N(A)>0,\ M(A)=0}}\hspace{25mm}\\
&\leq& \sum_{j=1}^n \P\pa{N^j(A)=1,\ \int_{A}\1_{X^j_{t-}=1}dt=0}\\
&\leq& \sum_{j=1}^n \P\pa{N^j(A)=1,\ \exists\varepsilon>0,\
  \lambda\pa{[T_{1,0}^j-\varepsilon,T_{1,0}^j]\cap A}=0}\\
&\leq& \sum_{j=1}^n \sum_{l\geq 1}\P\pa{T^j_{1,0}\in A,\
  \lambda\pa{[T_{1,0}^j-l^{-1},T_{0,1}^j]\cap A}=0}\\
&\leq& \sum_{j=1}^n \sum_{l\geq 1}\P\pa{T^j_{1,0}\in B_{l}}
\;\;=\;\;\sum_{j=1}^n \sum_{l\geq 1}\E\cro{N^*(B_l)}\;\;=\;\;0,
\end{eqnarray*}
by~\eref{eq:mark}. We may clearly fix $k=1$ and the choice of  $\kappa'$ is justified by the
following argument. First note that whatever $I\subset \X$ and
$t>0$
\begin{eqnarray*}
  \lefteqn{\P\pa{X^1_{t-}=1,\ T^1_{1,0}\in I,\ T^1_{1,0}\geq t}}\\
&=& \int_{I}\1_{\{u\ge t\}}\P\pa{X^1_{t-}=1,\ u\leq
  T^1_{1,0}\leq u+du}\\
&=& \int_{I}\1_{\{u\ge t\}}\P\pa{X^1_{t-}=1,\ X^1_{u-}=1}\P\pa{u\leq T^1_{1,0}\leq u+du
\,|\,X^1_{t-}=1,\ X^1_{u-}=1}\\
&=& \int_{I}\1_{\{u\ge t\}}\P\pa{X^1_{t-}=1,\ X^1_{u-}=1}\P\pa{u\leq T^1_{1,0}\leq u+du
\,|\,\ X^1_{u-}=1}
\end{eqnarray*}
since $X^1$ is a Markov process. Hence
\begin{eqnarray*}
 \P\pa{X^1_{t-}=1,\ T^1_{1,0}\in I,\ T^1_{1,0}\geq t}&=&\int_{I}
\1_{\{u\ge t\}}\P\pa{X^1_{t-}=1,\ X^1_{u-}=1}s(u)du\\ &=& \E\left[\int_I\1_{\{u\ge t\}}
\1_{\{X^1_{t-}=1\}}\1_{\{X^1_{u-}=1\}}s(u)du\right].
\end{eqnarray*}
It then follows that
\begin{eqnarray*}
{\rm Var}\pa{\int_I s(t)Y_t^1dt}&\le&\E\cro{\pa{\int_{\X}\1_I(t)s(t)Y_t^1dt}^2}\\&=&
\E\cro{\int_{\X\times\X}\1_{I}(t)\1_{I}(u)s(t)s(u)Y_t^1Y_{u}^1 dudt}\\&=& 2\int_I 
\E\cro{\int_{I}\1_{\{u\ge t\}}\1_{\{X^1_{t-}=1\}}\1_{\{X^1_{u-}=1\}}s(u)du}s(t)dt\\&=& 
2\int_I\P\pa{X^1_{t-}=1,\ T^1_{1,0}\in I, \ T^1_{1,0}\geq t}s(t)dt\\ &\leq&  2\int_I
\P\pa{X^1_{t-}=1}s(t)dt\;\;=\;\;2\E\cro{\int_I s(t)Y^1_tdt}.\\
\end{eqnarray*}

\mbox{ \vspace{10mm}}\\
\noindent\large{\bf {References\vspace{1mm}}}
\setlength{\parskip}{1mm}
{\small

ANDERSEN, P., BORGAN, O., GILL, R.\ and KEIDING, N. (1993). {\it Statistical Models
Based on Counting Processes}. Springer-Verlag, New York.

ANTONIADIS, A. (1989). A penalty method for nonparametric estimation of the
intensity function of a counting process. {\it  Ann. Inst. Statist. Math.} {\bf 41}, 781--807.

ANTONIADIS, A., BESBEAS, P.\ and SAPATINAS, T. (2001). Wavelet shrinkage for natural
exponential families with cubic variance functions. {\it Sankhya} {\bf 63}, 309-327.

ANTONIADIS, A.\ and SAPATINAS, T. (2001). Wavelet shrinkage for natural
exponential families with quadratic variance functions. {\it Biometrika} {\bf 88}, 805-820.

BARRON, A.R., BIRG\'E, L.\ and MASSART, P. (1999). Risk bounds for model selection via 
penalization. {\it Probab.\ Theory Relat.\ Fields} {\bf 113}, 301-415.

BARRON, A.R.\ and COVER, T.M. (1991). Minimum complexity density estimation. {\it 
IEEE Transactions on Information Theory}  {\bf37}, 1034-1054. 

BIRG\'E, L.  (1983). Approximation dans les espaces  m\'etriques   et th\'eorie de 
l'estimation. {\it Z. Wahrscheinlichkeitstheorie Verw. Geb.} {\bf65}, 181-237.  

BIRG\'E, L.  (2006). Model selection via testing : an alternative to (penalized) maximum
likelihood estimators. {\it Ann. Inst. Henri Poincar\'e Probab.\ et Statist.} {\bf 42}, 273-325.

BIRG\'E, L.\ and MASSART, P. (1998). Minimum contrast estimators on sieves: 
exponential  bounds and rates of convergence. {\it Bernoulli} {\bf 4}, 329-375.

BIRG\'E, L.\ and MASSART, P.  (2000). An adaptive compression algorithm in Besov
spaces.  {\it Constructive Approximation} {\bf 16} 1-36.

BIRG\'E, L.\ and MASSART, P. (2001). Gaussian model selection. {\it J. Eur. Math. Soc.}
{\bf 3}, 203-268.

BREIMAN, L., FRIEDMAN, J.H., OLSHEN, R.A.\ and STONE, C.J. (1984). {\it Classification
and Regression Trees}. Wadsworth, Belmont.

CASTELLAN, G. (1999). Modified Akaike's criterion for histogram density estimation. 
Technical Report 99.61. Universit\'e Paris-Sud, Orsay.

CASTELLAN, G. (2000). S\'election d'histogrammes \`a l'aide d'un crit\`ere de type
Akaike. {\it C.R.A.S.} {\bf 330}, 729-732.


DeVORE, R.A. (1998). Nonlinear Approximation. {\it Acta Numerica} {\bf 7}, 51-150.

DeVORE, R.A.\ and LORENTZ, G.G. (1993). {\it Constructive Approximation}. 
Springer-Verlag, Berlin. 

DeVORE, R.A.\ and YU,Ê X.M. (1990). Degree of adaptive approximation. {\it Math. Comp.}
{\bf  55}, 625-635.

GEY, S.\ and N\'ED\'ELEC, E. (2005). Model selection for CART regression trees.
{\it IEEE Transactions on Information Theory} {\bf  51}, 658-670.


GR\'EGOIRE, G.\and NEMB\'E, J. (2000). Convergence rates for the minimum
complexity estimator of counting process intensities. {\it
  J. Nonparametr. Statist.} {\bf  12}, 611-643.

KOLACZYK, E. (1999). Wavelet shrinkage estimation of certain Poisson intensity signals using
corrected threshold. {\it Statistica Sinica} {\bf 9}, 119-135.

KOLACZYK, E. \and NOWAK, R. (2004). Multiscale likelihood analysis and complexity
penalized estimation. {\it Annals of Statistics} {\bf 32}, 500-527.

LAURENT, B.\ and MASSART, P. (2000). Adaptive estimation of a quadratic functional 
by model selection. {\it Ann.\  Statist.}\ {\bf  28}, 1302-1338.

LEPSKII, O.V. (1991). Asymptotically minimax adaptive estimation I: Upper bounds.
Optimally adaptive estimates. {\it Theory Probab.\ Appl.}\  {\bf36}, 682-697.

MASSART, P.  (2000). Some applications of concentration inequalities to
Statistics. {\it Ann. Fac. Sciences de Toulouse}  {\bf IX},  245-303. 

PATIL, P.N.\ and WOOD, A.T. (2004). A counting process intensity estimation by
orthogonal wavelet methods.  {\it Bernoulli} {\bf 10}, 1-24.

REYNAUD-BOURET, P. (2003). Adaptive estimation of the intensity of
inhomogeneous Poisson processes via concentration inequalities. {\it Probab. Theory
Related Fields}  {\bf 126}, 103-153.

REYNAUD-BOURET, P. (2002). Penalized projection estimators of the Aalen
multiplicative intensity. {\it School of Mathematics, Georgia Institute of Technology}
Preprint 1202-002.

STANLEY, R.P (1999). {\it Enumerative Combinatorics, Vol.\ 2}. Cambridge University Press,
Cambridge.


van de GEER, S. (1995).  Exponential inequalities for martingales, with
application to maximum likelihood estimation for counting processes. {\it
Ann. Statist.} {\bf   23}, 1779-1801. 

WU, S.S.\ and WELLS, M.T. (2003) Nonparametric estimation of hazard functions by
wavelet methods. {\it  J. Nonparametr. Stat.} {\bf  15}, $187-203$.                                                                
}

\end{document}